\documentstyle{article}
 \input amssym.def
 \input amssym.tex

          \newtheorem{df}{Definition}[section]
          \newtheorem{pr}[df]{Proposition}
          \newtheorem{th}[df]{Theorem}
          \newtheorem{lem}[df]{Lemma}
          \newtheorem{cor}[df]{Corollary}
          \newtheorem{rem}[df]{\it Remark}

          \newtheorem{prob}[df]{Problem}
          \newtheorem{exam}[df]{Example}

          \newcommand{\qqed}{~$\Box$}
          
          \newcommand{\proof}{
                             \noindent{\bf Proof.~~}}
          \newcommand{\dslash}{\slash\hspace{-0.1cm}\slash}
          \newcommand{\va}{{\rm v}}
          \newcommand{\di}{{\rm div}}
            
            \makeatletter
            \@addtoreset{equation}{section}
            \makeatother
         
         \newcommand{\mapright}[1]{%
           \smash{\mathop{%
           \hbox to 1cm{\rightarrowfill}}\limits^{#1}}}
        
         \newcommand{\mapleft}[1]{%
           \smash{\mathop{%
           \hbox to 1cm{\leftarrowfill}}\limits_{#1}}}
           
         \newcommand{\mapdown}[1]{\Big\downarrow
         \llap{$\vcenter{\hbox{$\scriptstyle#1\,$}}$ }}

\begin{document}
      
         \title{Reduced class groups grafting relative invariants
        \thanks{Dedicated  to the memory of late Professor Masayoshi
        Nagata.}}
         
         \author{{ Haruhisa} {Nakajima}  \\
        \small Department of Mathematics, Faculty of Science, 
         {\sc  Josai University}\\
        \small Keyakidai, Sakado 350-0295, {\sc Japan}} 
         \date{ }
         
         \maketitle

         \begin{abstract}
            
             Let $(X, T)$ be a  regular stable conical 
              action of an algebraic torus on an affine normal conical variety $X$
              defined over an algebraically closed  field of characteristic zero. 
             We define  a certain  subgroup of   ${\rm Cl}(X\dslash T)$
               and 
              characterize its finiteness  in terms of
             a finite $T$-equivariant  Galois descent $\widetilde{X}$ of $X$. 
              Consequently we show that the  
              action $(X, T)$ is equidimensional if and only if
              there exists a 
              $T$-equivariant  finite Galois covering  $X\to \widetilde{X}$ such that
               $(\widetilde{X}, T)$ is cofree. Moreover
              the order of ${\rm Gal}(X/\widetilde{X})$ 
              is controlled by a certain subgroup of ${\rm Cl}(X)$. 
              The present  result
              extends thoroughly the equivalence  of  
              equidimensionality and   cofreeness  of $(X, T)$ for a factorial $X$. 
              The purpose of this paper is to evaluate orders of divisor classes associated
              to modules of relative invariants for a
              Krull domain with a group action.  This is useful  in 
              studying  on equidimensional torus actions as above.   
              The generalization of
              R. P. Stanley's criterion for freeness of modules of relative
              invariants plays an important role in showing  key assertions.
              
              
            \end{abstract}

          \small\tableofcontents
   \normalsize

         \small \section{Introduction}\normalsize
          Let $(X, G)$ denote a  regular action of an affine algebraic 
          group $G$ on an affine algebraic variety $X$ over an algebraically closed 
          field $K$ of characteristic $p\geq 0$.  We say that $(X, G)$ is  admitting
          an algebraic quotient, if the algebra ${\mathcal O}(X)^G$
           of invariants of $G$ in the coordinate ring ${\mathcal O}(X)$ of $X$ 
            is finitely generated over $K$. The algebraic quotient $X\dslash G$ 
            of $X$ by $G$
            is defined by ${\mathcal O}(X)^G$. 
          Such an action $(X, G)$ 
          is said to be {\it equidimensional} (resp. {\it cofree}), if  the quotient morphism
          $\pi_{X, G} : X \to X\dslash G$ 
          is equidimensional, i.e., closed fibers of  $\pi_{X, G}$ are
          pure ($\dim X -\dim X\dslash G$)-dimensional 
          (resp. 
          if ${\mathcal O}(X)$ is ${\mathcal O}(X)^G$-free). 
             
          For a finite dimensional linear representation $G \to GL(V)$
          of a connected  algebraic group $G$
          over the complex number field $\mbox{\boldmath $C$}$, the following conjecture
          is well known (e.g., Appendix to Chap. 4 of  \cite{Mum}) and  relates to
          the subject of the present paper (e.g., \cite{Nak2-1}).  
       
                 \medskip
          \noindent{\bf The Russian conjecture} {\it
          If $(V, G)$ is equidimensional, then it is cofree. }
          \medskip

           \noindent Especially for an  algebraic torus $G$, 
          the Russian  conjecture is solved affirmatively 
          (cf.  \cite{W}), which   is generalized 
          in  \cite{Nak1} to the case where $(V, G)$ is a conical factorial
          variety $V = X$ with a stable conical action of an algebraic torus
          over $K$ of characteristic zero. It should be
          noted that this result is not true for any conical normal variety
          $X$ (cf. Example 5.7, i.e., {\it the Russian conjecture 
          for normal varieties does not hold}\/). 
          Here $X$ (resp. $(X, G)$) is said to be {\it conical}, if $X$ is affine and
          ${\mathcal O}(X)$ is a positively graded,  i.e., a $\mbox{\boldmath $Z$}_0$-
          graded algebra defined over
          ${\mathcal O}(X)_0 = K$ 
          ($\mbox{\boldmath $Z$}_0$ = $\mbox{\boldmath $N$}\cup\{0\}$) (resp. if the action of $G$
          preserves each homogeneous part of ${\mathcal O}(X)$). Moreover
          $(X, G)$ is said to be {\it stable}, if there is a non-empty
          open subset of $X$ consisting of closed $G$-orbits. 
          
         The purpose of this paper is to study on the following  problem 
         which produces extensions of the results in \cite{Nak1}. 
          
          \begin{prob}\label{problem} Suppose that $G$ is an algebraic torus
          and $(X, G)$ a  stable conical action of $G$ on
          a conical normal variety $X$ defined over $K$
          of characteristic zero.  If  $(X, G)$
          is equidimensional, 
          then:
          \begin{itemize}
          \item{} Does  there  exist 
          a $G$-equivariant finite Galois covering $X\to \widetilde{X}$
          for a normal conical variety $\widetilde{X}$ with a
          conical $G$-action admitting the commutative diagram 
          \begin{eqnarray}
          X & \mapright{} & \widetilde{X}\nonumber\\
          \mapdown{\pi_{X, G}}& & \mapdown{\pi_{\widetilde{X}, G}}  \label{eqn1-1}\\
          X\dslash G & \mapright{\cong} & \widetilde{X}\dslash G \nonumber
         \end{eqnarray}  
          such that $(\widetilde{X}, G)$
          is   cofree?  
          \item{} Moreover can we choose $X\to \widetilde{X}$ in such a way that the order
           $\vert {\rm Gal} (X/\widetilde{X}) \vert$ of the Galois group of $X \to \widetilde{X}$ 
          is a divisor of  a power of the exponent  of a subgroup of the 
          divisor class group ${\rm Cl}(X)$ of $X$?  
          \end{itemize}
          \end{prob}

          This problem is closely connected with  relative invariants of algebraic tori, 
          because their irreducible representations are of degree one. 
          Relative invariants of finite groups and compact Lie groups are
          studied by R. P. Stanley (cf. \cite{St1, St2}), which inspires the author
          to  extend R. P. Stanley's criterion for freeness of  modules of relative invariants. 
          Our study on Problem \ref{problem} in the present paper
           is based on the result  in \cite{Nak0, Nak3} mentioned as above. 
           
           The main auxiliary part is given in   Sect. 3.  We explain the results
           in Sect.3 in more detail. Let us  consider
             a Krull domain $R$ acted by  an abstract group $G$  as ring automorphisms
             and   the $R^G$-modules $R_\chi$ of 
             $\chi$-invariants in $R$ for  $1$-cocycles $\chi$'s of $G$ in the
             group  ${\rm U}(R)$ of units of $R$. 
             We associate some qualified  cocycles $\chi$  with the Weil 
             divisors $D(\chi)$ on $R$ and 
             the divisorial $R^G$-lattices $d_{(R^G, R)}(R_\chi)$.
             In Corollary \ref{order of characteristic divisor}, the orders
             ${\rm ord} ([D(\chi)])$ in ${\rm Cl} (R)$ will be  characterized  by the equality   
             $${\rm ord} ([D(\chi)]) =
             {\rm min}\left(\left\{ m \in 
             \mbox{\boldmath $N$} \mid R_{m\chi} \cong R^G\right\}\right)$$
             on  the order of $[D(\chi)]$, 
             if  $R_{m(\chi)\chi} \cong R^G$ for some  $m(\chi)\in \mbox{\boldmath $N$}$.  
          By the identity  $$d_{(R^G, R)}(R_{m\chi}) = m\cdot
           d_{(R^G, R)}(R_{\chi})$$ for some cocycles, 
           we will establish  Theorem  \ref{main Krull case} as a main result in Sect. 3, 
           which   shows ${\rm ord} ([D(\chi)])$ is equal to
           ${\rm ord} ([d_{(R^G, R)}(R_{\chi})])$ in ${\rm Cl}(R^G)$. In Definition 
           \ref{theta}, {\it the reduced
           class group} ${\rm UrCl}(R, G)$ (resp.  $\widetilde{{\rm Cl}}(R, G)$) 
           is  defined to be a  subgroup of ${\rm Cl}(R)$ (resp. of ${\rm Cl}(R^G)$)
           generated by certain $[D(\chi)]$'s (resp.  certain 
           $[d_{(R^G, R)}(R_\chi)]$'s).  Then
           by Theorem \ref{finiteness of class group} we have the equality of 
           exponents 
           \begin{equation}
           {\rm exp} ({\rm UrCl}(R, G)) = {\rm exp} (\widetilde{{\rm Cl}}(R, G)) <\infty
           \label{eqn1-2}\end{equation}
           of reduced class groups 
           if $R_{m(\chi)\chi} \cong R^G$ for some qualified $\chi$'s with $m(\chi) \in 
           \mbox{\boldmath $N$}$.

           Sect. 4, 5  are devoted to the study on Problem \ref{problem} and hence,  without
           specifying,  suppose 
            that $(X, G)$ is a regular action of a connected
           affine algebraic group $G$ on an affine normal variety over $K$ of characteristic zero. 
          The  pseudo-reflections in the general linear group
           $GL(V)$ are recognized as elements 
           which are inertial  at minimal prime ideals in ${\mathcal O}(V) =
           {\rm Sym} (V^\ast)$ (cf. \cite{Nak0}). We define pseudo-reflections
           for $(X, G)$ and 
           characterize the ramification indices of discrete valuations 
           on $X$ over $X\dslash G$ in the case where $G^0$  is
           an algebraic torus,  
           in terms of orders of pseudo-reflections (cf. \cite{Nak2}). The qualified   
           cocycles treated in Sect. 3 are determined by pseudo-reflection
           subgroups (cf.
           Proposition \ref{pseudo-reflection group}). Moreover
           suppose that both $X$ and $(X, G)$ are conical. Then, combining  
           Proposition \ref{pseudo-reflection group} 
           with Theorem \ref{finiteness of class group},  by \cite{Nak1}
           we will obtain
           the main result of Sect. 4, 5,   i.e., Theorem \ref{main theorem for algebraic groups}.
           Furthermore we define in Definition \ref{obstruction}
            the  {\it obstruction subgroup} ${\rm Obs}(X, G)$ of $G$  {\it for cofreeness of} 
            $(X, G)$ 
           which excludes some characters $\chi$ such that ${\mathcal O}(X)_\chi\not\cong
           {\mathcal O}(X)^G$.  Consequently, 
           we solve affirmatively Problem \ref{problem}  as follows (cf. Theorem \ref{equi}):

          \begin{th}\label{1-2} Under the same circumstances as in Problem \ref{problem},   
          $(X, G)$ is equidimensional if and only if there exists a  
           (normal) closed subgroup ${\rm Obs}(X, G)$ of $G$ whose restriction
           to $X$ is of   order dividing 
          a power of the exponent of the 
             subgroup ${\rm UrCl}({\mathcal O}(X), G)$ of ${\rm Cl}(X)$
          such that $(X\dslash {\rm Obs}(X, G), G)$ is cofree. 
          \end{th}
          
         The order
         $\vert {\rm Obs}(X, G)\vert_X\vert$ of restrictions of ${\rm Obs}(X, G)$ to
         $X$ is effectively determined and is closely related to 
          the equality (\ref{eqn1-2}). 
         Clearly putting $\widetilde{X} : = X\dslash {\rm Obs}(X, G)$, we
         have a commutative diagram (\ref{eqn1-1}) and 
         ${\rm Gal}(X\slash \widetilde{X}) = {\rm Obs}(X, G)\vert_X$.  
         Thus  the following criterion for cofreeness 
        is  obtained:

          \begin{cor}\label{1-3}
          Under the same circumstances as in Problem \ref{problem},
          suppose that $(X, G)$ is equidimensional. Then
          $(X, G)$ is cofree if and only if $ \left\vert{\rm Obs}(X, G)\vert_X\right\vert
          =1 $.  
           \end{cor}
         The {\it if part} of Corollary \ref{1-3} is an immediate consequence of Theorem \ref{1-2}.  
         There are many equidimensional $(X, G)$'s which are not cofree (e.g., 
         Example \ref{example2}). 
          Theorem \ref{1-2} and Corollary \ref{1-3} are regarded as
           quite generalizations of the main theorem
         in \cite{Nak1}. 
         
              For a group homomorphism $G \to {\rm Aut}(Y)$ 
        we
          denote by $G\vert_ Y$ the set $\{ \sigma\vert_Y \mid \sigma \in G \}$
          where $\sigma\vert_Y$ denotes the restriction of $\sigma$ to $Y$. 
          Let ${\rm exp} (G)$ be the {\it exponent} of a group $G$
          and ${\rm ord} (a)$ denote the {\it order} of $a\in G$. For a subset
         (or an element)  $\Omega$ of $G$, let $\left< \Omega \right>$ be 
          the  subgroup of $G$ generated by $\Omega$ and $Z_G(\Omega)$ the
          centralizer of $\Omega$ in $G$. 
          Let ${\rm tor}(A)$ denote the {\it torsion part} of an abelian group $A$.
           The notations $\mbox{\boldmath $Z$}$
          and $\mbox{\boldmath $N$}$ are  standard (i.e, the
          set of all integers and that of all natural numbers, respectively). 
           Let $\mbox{\boldmath $Z$}_0$
          denote the additive monoid consisting of all non-negative integers.

         \small \section{Preliminaries} \normalsize

           Let ${\mathcal Q}(A)$ be the total quotient ring of a 
          commutative ring $A$ and let ${\rm Ht}_1(A)$ be the set consisting of 
          all prime ideals ${\frak p}$ of $A$ of height one (${\rm ht} ({\frak p}) = 1$).
          Consider  a ring extension $A \hookrightarrow B$ of 
          integral domains.  For $i\in \mbox{
          \boldmath $Z$}_0$ set
          \[{\rm Ht}_1^{(i)}(B, A) := \{{\frak p}\in {\rm Ht}_1(B) \mid {\rm ht} ~({\frak p}\cap A)
          = i \}\] and moreover 
          \({\rm Ht}_1^{[j]}(B, A) := \bigcup _{i\geqq j} {\rm Ht}_1^{(i)} (B, A)\)
          for $j\in \mbox{\boldmath $N$}$. 
          Put  $$X_{\frak Q}(B):=
           \{{\frak p}\in {\rm Ht}_1(B) \mid {\frak p}\cap A={\frak Q}\}$$
           for any ${\frak Q}\in{\rm Ht}_1(A)\). 
          For a nonempty subset  $M$  of  ${\mathcal Q}(B)$, we define   $${ d}_{(A, B)}(M):=
          \bigcap _{{\frak q}\in {\rm Ht}_1(A)} \left(A_{\frak q}\otimes _A (A\cdot M)
          \right).$$  Here $A\cdot M$ denotes the $A$-submodule
          generated by $M$ and  each  $A_{\frak q}\otimes _A (A\cdot M)$ 
          is regarded as a subset of ${\mathcal Q}(A)\otimes _A (A\cdot M)
          \hookrightarrow {\mathcal Q}(B)$. Especially in the case of $A =B$, 
          we simply use $d_A(M)$ instead of the notation $d_{(A, A)}(M)$. 
           If $B$ is a Krull domain
          and $A = {\mathcal Q}(A)\cap B$, then $A$ is a Krull domain
          (e.g., \cite{Matsumura}) and for any ${\frak q}\in {\rm Ht}_1(A)$
          the set $X_{\frak q}(B)$ is nonempty (cf. \cite{Magid, Nak1}). In this case, for
          ${\frak P}\in X_{\frak q}(B)$ let ${\rm e}({\frak P}, {\frak q})$
          denote the reduced ramification index of ${\frak P}$ over ${\frak q}$
          (e.g., \cite{Nak2}). 
          In the case where $A$ is a Krull domain, let $\va_{A, {\frak p}}$ denote
          the discrete valuation of $A$ defined by ${\frak p}\in {\rm Ht}_1(A)$ and,
          for a subset $M$ of ${\mathcal Q}(A)$ generating a fractional
          ideal $A\cdot M$ of $A$ in ${\mathcal Q}(A)$, let $\di_A(M)$ be
          the divisor on $A$ associated with the divisorialization $d_A(A\cdot M)=
          d_A(M)$
          of $A\cdot M$ on $A$ (e.g., \cite{Bass}) . Let ${\rm Div}(A)$   
          and ${\rm Cl}(A)$ denote
          the (Weil) divisor group
          and the divisor class
          group of $A$, respectively. For a (Weil) divisor $D \in {\rm Div}(A)$, 
          we denote by $I_A(D)$ the divisorial fractional ideal
          of $A$ associated with $D$. Our notation is standard (cf. \cite{Bass, Fossum}).

          \begin{pr}[\cite{Nak3}]\label{n1} Let $R$ be a Krull domain 
          and $L$ a subfield of ${\mathcal Q}(R)$ such that $L = {\mathcal Q}(R\cap L)$.
          Suppose that ${\frak J}$ is a divisorial fractional ideal of $R$. If  ${\frak J}\cap L
          \not= {0}$, then ${\frak J}\cap L$ is a fractional ideal of $R\cap L$
          and $\di_{R\cap L}( {\frak J}\cap L)$ equals to $D_{\frak J}$ which
          is given by 
          \[D_{\frak J}:=
          \sum_{{\frak q}\in {\rm Ht}_1(R\cap L)} \Bigl( \max_{{\frak P}\in X_{\frak q}(R)}
          \Bigl[\frac{\va_{R, {\frak P}}({\frak J})}{{\rm e}({\frak P}, {\frak q})} 
          \Bigr]^\sharp\Bigr)\cdot \di_{R\cap L}({\frak q}), \]
          where $[a]^\sharp$ denotes 
          $-[-a]$ for any real number $a$
          and $[~\cdot~]$ denotes the Gauss symbol. 
          
          \end{pr}
          
          \begin{cor}\label{power of ideal}
           Under the same circumstances  as in Proposition \ref{n1},
          suppose moreover that \({\rm v}_{R, {\frak P}} ({\frak J}) \equiv 0~
          {\rm mod} ~{\rm e} ({\frak P}, {\frak q})\)
          for any ${\frak q}\in {\rm Ht}_1(R\cap L)$ and any ${\frak P} \in X_{\frak q}(R)$.
          Then 
          \[ n \cdot \di_{R\cap L}( {\frak J}\cap L) = \di_{R\cap L}
          (d_R({\frak J}^n)\cap L)\]
          for any $n \in \mbox{\boldmath $N$}$.
          \end{cor}
          
          \proof Let  ${\frak q} \in {\rm Ht}_1 (R\cap L)$ be any prime ideal. 
          For  a prime ideal ${\frak P_0}\in X_{\frak q}(R)$,   the
          condition 
          \[\max _{{\frak P}\in X_{\frak q}(R)} ({\rm v}_{R, {\frak P}} ({\frak J}))
          = {\rm v}_{R, {\frak P_0}} ({\frak J})\] 
          holds if and only if 
          \[\max _{{\frak P}\in X_{\frak q}(R)} ({\rm v}_{R, {\frak P}} (d_R({\frak J}^n)))
          = {\rm v}_{R, {\frak P_0}} (d_R({\frak J}^n)).\] 
          By the assumption on discrete 
          valuations of ${\frak J}$ in $R$ and Proposition \ref{n1},
          we have
          \begin{eqnarray} {\rm v}_{R\cap L, {\frak q}}(d_R({\frak J}^n)\cap L)
          &=& \max _{{\frak P}\in X_{\frak q}(R)} \Bigl(
          \frac{{\rm v}_{R, {\frak P}} (d_R({\frak J}^n))}{{\rm e} ({\frak P}, {\frak q})}\Bigr)
          =\frac{{\rm v}_{R, {\frak P_0}} (d_R({\frak J}^n))}
          {{\rm e} ({\frak P}, {\frak q})}\nonumber\\
          &=& n\cdot  \max _{{\frak P}\in X_{\frak q}(R)} \Bigl(
          \frac{{\rm v}_{R, {\frak P}} ({\frak J})}{{\rm e} ({\frak P}, {\frak q})}
          \Bigr)\nonumber\\
          &=& n\cdot {\rm v}_{R\cap L, {\frak q}}({\frak J}\cap L).\nonumber
          \end{eqnarray}
           The assertion follows from this. \qqed
          
          \medskip
               
      For a Krull domain $R$ and a subfield $L$ of ${\mathcal Q}(R)$,  
      the subgroup $D_{(R, R\cap L)}(g_1, \dots, g_n)$ of ${\rm Div}(R\cap L)$ 
      is defined as follows:
      
      \begin{lem} \label{finite lemma}  Suppose that $R$ is a Krull domain and
      let $L$ be a subfield of ${\mathcal Q}(R)$ such that $L = {\mathcal Q}(R\cap L)$.
     Let    $\{ { g}_1, \dots {g}_n\}$ be a finite
     set  of nonzero elements   in $R$.  Then the
           subgroup 
           \[ \left< \left\{\di _{R\cap L} \Bigl(\Bigl(\left. \frac{1}{\prod _{j=1}^n 
           {g_{j}}^{i_j} } R\Bigr) \cap L
           \Bigr) ~\right\vert~ \forall i_1, \dots, i_n \in \mbox{\boldmath $Z$}_0 \right\} \right>,
           \]
           which is denoted by $D_{(R, R\cap L)}(g_1, \dots, g_n)$,
            of ${\rm Div} (R\cap L)$ is finitely generated.

          \end{lem}
          
          \proof By Proposition \ref{n1}, we
          see that, for any $\frak{q}\in {\rm Ht}_1(R\cap L)$, the equivalence
          \[{\rm v}_{R\cap L, \frak{q}}\Bigl(\Bigl(\frac{1}{\prod _{j=1}^n 
           {g_{j}}^{i_j} } R\Bigr) \cap L  \Bigr) \not= 0 
                      \Longleftrightarrow
           \max_{{\frak P}\in X_{\frak q}(R)} 
          \left[-~\frac{\sum_{j=1}^n i_j \cdot \va_{R, {\frak P}}
          ({R\cdot {g_{j}}})}{{\rm e}({\frak P}, {\frak q})} 
          \right]^\sharp \not= 0
          \]
           holds,  where $[~ \cdot~]^\sharp$ defined in Proposition \ref{n1}. 
                   On the other hand,  the set
          \[\bigcup _{j=1}^n \Bigl\{{\frak q}\in {\rm Ht}_1(R\cap L) ~\Bigm
          \vert~ \exists {\frak P}\in 
          X_{\frak q}(R)~
          \mbox{such that}~ {\rm v}_{R, {\frak P}}( R\cdot g_j)\not= 0\Bigr\},\]
         which is denoted by ${\rm supp} (g_1, \dots, g_n: R, R\cap L)$, 
         has only finite elements.  Since
         \[\left\{ {\frak q} \in \in {\rm Ht}_1(R\cap L) ~\left\vert~
         {\rm v}_{R\cap L, \frak{q}}\Bigl(\Bigl(\frac{1}{\prod _{j=1}^n 
           {g_{j}}^{i_j} } R\Bigr) \cap L  \Bigr) \not= 0\right. \right\}
           \subseteq {\rm supp} (g_1, \dots, g_n: R, R\cap L)
           \]
           for any $i_1$, $\dots$, $i_n\in \mbox{\boldmath $Z$}_0$, we see that
           \[D_{(R, R\cap L)}(g_1, \dots, g_n) \subseteq \sum_{{\frak q}\in {\rm supp} 
           (g_1, \dots, g_n: R, R\cap L)}
           \mbox{\boldmath $Z$} \cdot {\rm div}_{R\cap L} ({\frak q})
           \subseteq {\rm Div} (R\cap L).\]
           The assertion follows from this. \qqed
           
           \medskip

          Suppose that $R$ is  an integral  domain   
          on which a group $G$ acts as automorphisms. 
          Let $Z^1(G, {\rm U}(R))$ be the group 
          of all $1$-cocycles of $G$ on the group ${\rm U}(R)$ of units of $R$ whose
          group structure  is given by addition.  The {\it trivial $1$-cocyle} is 
          denoted by $\theta\in Z^1(G, {\rm U}(R))$.
          For any $\chi\in Z^1(G, {\rm U}(R))$,
          let $$R_\chi = \{ x\in R \mid \sigma (x) = \chi (\sigma) \cdot x
          ~~(\forall \sigma \in G)\}$$
          which is regarded as an $R^G$-module. If $\Lambda$ and $\Gamma$
          are subsets of $Z^1(G, {\rm U}(R))$, let   $-\Lambda$
           be the set $\{-\chi  \in Z^1(G, {\rm U}(R)) \mid \chi \in \Lambda\}$
         and put $$\Lambda + \Gamma := \{\chi +\psi \in Z^1(G, {\rm U}(R)) 
         \mid \chi \in \Lambda, 
         \psi \in \Gamma\}.$$    
          Denote by  $Z^1(G, {\rm U}(R))^R$ the
          set $\{\chi \in Z^1(G, {\rm U}(R)) \mid R_\chi \not= \{0\} \}$ and put  
           \begin{eqnarray*} & Z^1_R(G, {\rm U}(R))_{(2)}:=
          \Bigl\{ \chi \in - Z^1(G, {\rm U}(R))^R~\Bigm\vert ~   
           R_{-{\chi}} \not\subseteq {\frak P}~ ~ (\forall
          {\frak P}\in{\rm Ht}_1^{[2]}(R, R^G) ) \Bigr\}, \end{eqnarray*}
          $$\widetilde{Z}^1_R(G, {\rm U}(R)):= Z^1_R(G, {\rm U}(R))_{(2)}\cap
          \left(- Z^1_R(G, {\rm U}(R))_{(2)}\right)$$ 
          respectively.  Set 
          $$Z^1(G, {\rm U}(R))_R := \Bigl\{ 
          \chi \in Z^1(G, {\rm U}(R))^R ~\Bigm\vert~ \dim_{{\mathcal Q}(R^G)} {\mathcal Q}(R^G)
          \otimes_{R^G} R_\chi = 1\Bigr\}.$$  
          For $\chi \in Z^1(G, {\rm U}(R))^R$, the condition 
          that 
          $\dim_{{\mathcal Q}(R^G)} {\mathcal Q}(R^G)
          \otimes_{R^G} R_\chi = 1$ holds if and only if 
          the equality
          $$((1/f) \cdot R) \cap {\mathcal Q}
          (R^G ) = ((1/f) \cdot R)^G$$ holds for some (or any) nonzero $f \in R_\chi$ 
         (cf. Lemma 3.1 of \cite{Nak3}). 
          Moreover in  the case
          where $R$ is a Krull domain, put 
           \begin{eqnarray*}  Z^1_R(G, {\rm U}(R))_{\rm e}:=
           \Bigl\{ \chi\in Z^1(G, {\rm U}(R))^R ~\Bigm\vert~ \exists f_{\frak P}
          \in R_\chi\backslash \{0\} ~\mbox{such that}~\\
           {\rm v}_{R, {\frak P}} (f_{\frak P}) 
         \equiv 0~ ({\rm mod} ~{\rm e} ({\frak P}, {\frak P}\cap R^G))~~
          (\forall {\frak P}\in {\rm Ht}_1^{(1)}(R, R^G)) \Bigr\}. \end{eqnarray*}

       \begin{lem}\label{easy property} 
       The set $Z^1(G, {\rm U}(R))^R$ has the following properties:
       \begin{itemize}
       \item[{\it (i)}]
       For $\chi_i \in Z^1(G, {\rm U}(R))^R$ ($i =1, 2$), the condition 
       $\chi_1 + \chi_2 \in Z^1(G, {\rm U}(R))_R$
       implies that $\chi_i \in Z^1(G, {\rm U}(R))_R$. 
       \item[{\it (ii)}]
        $Z^1(G, {\rm U}(R))_R \supseteq Z^1(G, {\rm U}(R))^R\cap
          \left(- Z^1(G, {\rm U}(R))^R\right)$.  
          \item[{\it (iii)}] 
          If  ${\mathcal Q}(R^G) = {\mathcal Q}(R)^G$, 
          then $Z^1(G, {\rm U}(R))_R = Z^1(G, {\rm U}(R))^R$. 
          \end{itemize}
          
          \end{lem}
          \proof ~{\it (i):} Let $f$ be a nonzero element of $R_{\chi_2}$. Then 
          the map $R_{\chi_1} \ni x \to x\cdot f \in R_{\chi_1+\chi_2}$ is
          $R^G$-monomorphism, which implies the assertion. 
          \par The remainder of the assertions follows from Lemma 3.1 of \cite{Nak3}. 
          \qqed
       
       \medskip
          
          We immediately have  the following lemma  whose proof is easy and 
          omitted:

          \begin{lem}\label{semigroup} The set
          $Z^1_R(G, {\rm U}(R))_{(2)}$ (resp. $\widetilde{Z}^1_R(G, {\rm U}(R))$) 
          is an additive submonoid (resp. subgroup) of $Z^1(G, {\rm U}(R))$ 
          containing the group $B^1(G, {\rm U}(R))$ of  
          $1$-st coboudaries of $G$ in ${\rm U}(R)$.
          Especially in the case where $R$ is a Krull domain,  the set
          $$Z^1_R(G, {\rm U}(R))_{(2)}\cap Z^1_R(G, {\rm U}(R))_{\rm e} ~\mbox{ 
          (resp. }~\widetilde{Z}^1_R(G, {\rm U}(R))\cap Z^1_R(G, {\rm U}(R))_{\rm e})$$ 
          is a submonoid (resp. subgroup) of $Z^1(G, {\rm U}(R))$ 
          containing the group $B^1(G, {\rm U}(R))$.
          \qqed
         \end{lem}

          \begin{lem}\label{divisorialization}  
           Let $\chi$ be a 1-cocycle
          in  $Z^1(G, {\rm U}(R))^R$. 
          \begin{itemize}
          \item[(i)] 
          $R_\chi ={d}_{(R^G, R)} (R_\chi)$ if and only if 
          ${d}_{(R^G, R)} (R_\chi)\subseteq R$.
           \item[(ii)] Suppose that  $\chi \in Z^1(G, {\rm U}(R))_R$  
           and  ${d}_{(R^G, R)} (R_\chi)\subseteq R$ hold. Then 
          \(R_\chi\cong R^G\) as \(R^G\)-modules  if and only if 
          the following condition is satisfied:
           For  a nonzero element $f$ of $R_\chi$, the $R^G$-module 
           ${d}_{(R^G, R)}
           \displaystyle\left( \left((1\slash f) \cdot R\right) \cap {\mathcal Q}(R^G)\right)$
           is  principal.

          \end{itemize}
          \end{lem}

                     \proof ~{\it (i):}~ The {\it``only if"} part
                     follows from the fact 
                     $d_{(R^G, R)}(R_\chi) \subseteq ({\mathcal Q} (R^G)
                     \otimes_{R^G} R)_\chi$. 
                     
                     {\it (ii):}~By the choice of $\chi$ and
                     Lemma 3.1 of \cite{Nak3}, we see 
                     $${d}_{(R^G, R)}
           \displaystyle\left( \left((1\slash f) 
           \cdot R\right) \cap {\mathcal Q}(R^G)\right) =
           {d}_{(R^G, R)}
           \displaystyle\left((1\slash f) 
           \cdot R_\chi\right) =  (1\slash f) 
           \cdot  {d}_{(R^G, R)}
           \displaystyle\left(R_\chi\right),$$
           which shows the assertion. \qqed


          \small\section{Krull domains with group actions} \normalsize

          In this section suppose that $R$ is a Krull domain
          acted by  a group $G$  as automorphisms. Let $E^\ast (G, R)$ denote
          the subgroup $$\displaystyle \sum_{{\frak q}\in {\rm Ht}_1(R^G)}
          \mbox{\boldmath $Z$} \cdot \Bigl(\sum_{{\frak P}\in X_{\frak q}(R)} 
          {\rm e} ({\frak P}, {\frak q}) \cdot {\rm div}_{R}({\frak P})
          \Bigr) \oplus \Bigl( \sum _{{\frak P}\in {\rm Ht}_1^{[2]}(R, R^G)
          } \mbox{\boldmath $Z$} \cdot {\rm div}_{R}({\frak P})\Bigr)$$ 
          of ${\rm Div} (R)$. Recall that a divisor $D$ is said to be 
          effective (i.e.,  $D\geqq 0$), if $I_R(D) \subseteq R$. 
          \begin{df} \rm An effective
          divisor $D\in {\rm Div} (R)$, which is denoted by  $D\geqq 0$, 
            is said to be {\it minimal effective relative to} $(R^G, R)$,
          if $D$ has a decomposition $D = D_1 + D_2$ for $0 \leqq D_1\in E^\ast (G, R)$ 
          and $0 \leqq D_2\in
          {\rm Div} (R)$, 
          then  $D_1$ must  be equal to  zero.  
          \end{df}
          
          For a minimal
          effective divisor $D$ relative to $(R^G, R)$ on $R$, we immediately see 
          $D = 0$ if and only if $D \equiv 0 ~{\rm mod}~ E^\ast(G, R)$.  
          In the following lemma
          {\it we define the minimal effective divisor $D(\chi)$  relative to
          $(R^G, R)$ for any cocycle $\chi \in Z^1(G, {\rm U}(R))_R$.} 
          The group $G$ acts on ${\rm Div} (R)$ naturally and a divisor $D$ on $R$
          is said to be $G$-invariant  if $D$ is invariant under the action of $G$ 
          on ${\rm Div} (R)$,  i.e., $D\in {\rm Div} (R)^G$.

          \begin{lem}
          \label{characteristic divisor} 
       Let  $\chi$ be a cocycle
          in $Z^1(G, {\rm U}(R))_R$.  Then: 
           \begin{itemize}
           \item[(i)] For any nonzero $f$, $g\in R_\chi$, 
           ${\rm div}_R(f) \equiv {\rm div}_R(g) ~ {\rm mod}~E^\ast (G, R).$
            \item[(ii)]  If  $\chi\in 
          Z^1_R(G, {\rm U}(R))_{\rm e}$, then for any nonzero $f\in R_\chi$ we have
          $${\rm v}_{R, {\frak P}} (f) \equiv 0~ ({\rm mod}~ 
          {\rm e} ({\frak P}, {\frak P}\cap {\frak q}))~
          (\forall {\frak P}\in {\rm Ht}_1^{(1)}(R, R^G)).$$
          \item[(iii)]
          There exists a unique minimal effective divisor $D(\chi)$ on $R$
          relative to $(R^G, R)$ 
          such that, for a nonzero element $f \in R_\chi$, 
          $$E^\ast (G, R) \ni {\rm div}_R (f) -D(\chi) \geqq 0.$$ 
          \item[(iv)]  The divisor $D(\chi)$ 
          is   $G$-invariant  
          and  does not depend on the choice of   $f\in R_\chi$.  
          Moreover the   divisorial ideal $I_R(D(\chi))$ associated with
          $D(\chi)$ contains 
          $R_\chi$.  
         \item[(v)]  In the case where  $\chi \in Z^1_R(G, {\rm U}(R))_{\rm e}$, we have 
         $${\rm v}_{R, \frak P}(I_R(D(\chi))) \equiv 0 ~({\rm mod}~ 
          {\rm e} ({\frak P}, {\frak P}\cap {\frak q}))~
          (\forall {\frak P}\in {\rm Ht}_1^{(1)}(R, R^G)).$$ Moreover if
         $m\chi \in Z^1(G, {\rm U}(R))_R$ for some 
         $m \in \mbox{\boldmath $N$}$, then 
           the divisors $D(\chi)$ and $D(m \chi)$
          satisfy 
          $m\cdot D(\chi) = D(m \chi)$ in ${\rm Div} (R)$.
          \item[(vi)] For any cocycles $\chi_1$, $\chi_2 \in Z^1(G, {\rm U}(R))^R$
          such that $\chi_1 + \chi_2\in Z^1(G, {\rm U}(R))_R$,
         the divisors $D(\chi_1)$ and $D(\chi_2)$ (for existence of $D(\chi_i)$, 
          see {\it (i)} of 
          Lemma \ref{easy property})  
          satisfy 
          $$D(\chi_1 + \chi_2) \equiv D(\chi_1) +D(\chi_2) ~{\rm mod}~E^\ast (G, R).$$
          \item[(vii)] Suppose that  $D(\chi)$ is principal, i.e., 
          $D(\chi) =  {\rm div}_R (g)$ for some $g \in R$. 
           Then there is  a
         cocycle $$\psi\in \bigr(\chi + (-Z^1(G, {\rm U}(R))^R)
         \bigl)\cap Z^1(G, {\rm U}(R))^R$$
         such that  $g\in R_\psi$. This cocycle  $\psi$ has the properties 
         that $\psi \in Z^1(G, {\rm U}(R))_R$, $D(\psi) = D(\chi)$ and $D(\chi -\psi) = 0$.

          \end{itemize}

        \end{lem} 
        \proof  Since the choice of $\chi$ 
        implies $\displaystyle {f}\slash {g} = {{f}'}\slash {g'}$ for
        some nonzero $f'$, $g' \in R^G$, 
        the assertion {\it (i)} follows easily from this.
        The assertion {\it (ii) } is a consequence of {\it (i) } and the definition of 
        $ Z^1_R(G, {\rm U}(R))_{\rm e}$. 
        
        {\it (iii):} ~For any nonzero $f \in R_\chi$, put
         \begin{eqnarray}
        D_2(f)  
         &=& \sum_{{\frak P}\in {\rm Ht}_1^{(0)}(R, R^G)
         \cup {\rm Ht}_1^{(1)}(R, R^G) }{\rm v}_{R, {\frak P}}(f)\cdot
         {\rm div}_R({\frak P}) \nonumber\\
         &-& \sum_{{\frak q}\in {\rm Ht}_1(R^G)} \Bigl(
         \min_{{\frak Q}\in X_{\frak q}(R)}\left[ \frac{{\rm v}_{R, {\frak Q}}(f)}
         {{\rm e}({\frak Q}, {\frak q})}\right]\Bigr)\cdot
         \Bigl(\sum_{{\frak P}\in X_{\frak q}(R)} {\rm e}({\frak P}, {\frak q})\cdot
         {\rm div}_R({\frak P})\Bigr)\label{D2}
         \end{eqnarray}
         and $D_1(f) = {\rm div}_R(f) - D_2(f)$, where $[~\cdot~]$ denotes
         the  Gauss symbol. Then we see $0\leqq D_1(f) \in 
         E^\ast (G, R)$ and 
         $D_2(f)$ is a minimal  
         effective divisor on $R$ relative to $(R^G, R)$. 
         For a nonzero $g \in R_\chi$,  define $D_1(g)$ and $D_2(g)$ similarly as above.
          By {\it (i)} we see that 
         \begin{eqnarray} D_2(f) - D_2 (g) &\equiv&  
         - D_1(f) + D_1(g) ~{\rm mod}~ E^\ast (G, R)\nonumber \\
         &\equiv& 0  ~{\rm mod}~ E^\ast (G, R)\label{D3}\end{eqnarray}
         which implies, for any  ${\frak q}\in {\rm Ht}_1(R^G)$,  there
         is an integer $c_{\frak q}$ satisfying
         $${\rm v}_{R, {\frak P}}(I_R(D_2(f))) - {\rm v}_{R, {\frak P}}(I_R(D_2(g))) = c_{\frak q}
         \cdot {\rm e}({\frak P}, 
         {\frak q})\hskip 0.3cm (\forall {\frak P}\in X_{\frak q}(R))$$
         and depending only on ${\frak q}$. 
         Exchanging $f$ with $g$, we may suppose that $c_{\frak q} \geqq 0$.
         By the definition (\ref{D2}) of $D_2(f)$, we see 
         ${\rm v}_{R, {\frak P}}(I_R(D_2(f))) < {\rm e}({\frak P}, 
         {\frak q})$ for some ${\frak P}\in X_{\frak q}(R)$, which requires
         $c_{\frak q} = 0$. Consequently  $D_2(f)$ is 
         just the  minimal effective divisor
         desired in the assertion {\it  (iii)}. 
         
         {\it (iv)}: By  the congruence (\ref{D3}) the divisor $D_2(f)$ 
           does not depend on the choice of a nonzero
         $f\in R_\chi$. The remainder  follows from {\it (iii)}. 
         
         {\it (v):}~The congruence in {\it (v)}
          follows from {\it (ii)} and 
         (\ref{D2}). 
         Choose any  ${\frak q}$ from ${\rm Ht}_1(R^G)$. 
         Then, by the choice of $\chi$, the assertion {\it  (ii)} and the 
         definition of $D(\chi)$,  we see that ${\rm v}_{R, {\frak P}} (I_R(D(\chi)))
         = 0$ for some ${\frak P}\in X_{\frak q}(R)$, which implies 
         ${\rm v}_{R, {\frak P}} (I_R(m\cdot D(\chi)))
         = 0$.  
          From this 
         we see that $m\cdot D(\chi)$
         is a minimal effective divisor on $R$ relative to
         $(R^G, R)$. Since $$E^\ast (G, R) \ni m \cdot \left({\rm div}_R(f)- 
         D(\chi) \right) = 
         {\rm div}_R(f^m) - m\cdot D(\chi) \geqq 0,$$
          the assertion follows from
         the uniqueness of $D(m\chi)$ for $m\chi \in  Z^1 (G, {\rm U}(R))_R$
          shown  in the assertion {\it (iv)}. 
          
         {\it (vi):}~ Note that 
        $\chi_i \in Z^1(G, {\rm U}(R))_R$ (cf. {\it (i)} of
        Lemma \ref{easy property}) and $D(\chi_i)$ are
        well defined. The congruence in {\it (vi)} 
        follows easily from {\it (iii)} and the independence of $D(\chi)$ on the choice
        of $f$ stated in {\it (iv)}. 
    
       {\it (vii):}~ 
       As $D(\chi)$ is $G$-invariant, we have a cocycle $\psi\in Z^1(G, {\rm U}(R))^R$
       such that $R_\psi \ni g$. Then $D(\chi)$ is a minimal effective 
       divisor relative to $(R^G, R)$ with the property
       $$E^\ast (G, R) \ni {\rm div}_R(g) - D(\chi) = 0 ~( \geq 0).$$
        Since $\{0 \}\not=
       R_\chi \subseteq I_R(D(\chi)) = R\cdot g$
       (cf. {\it (iv)}), we see that $\chi -\psi \in Z^1(G, {\rm U}(R))^R$, which
       implies $\psi$, $\chi -\psi \in Z^1(G, {\rm U}(R))_R$ (cf. {\it (vi)}). 
       Clearly  by the uniqueness of $D(\psi)$ (cf. {\it (iii)}), 
       $D(\chi) = D(\psi)$.
       By {\it (vi)} we have $$D(\psi) + D(\chi -\psi) \equiv D(\chi) ~{\rm mod}
       ~E^\ast(G, R).$$
         This shows  $D(\chi - \psi) = 0$.    \qqed
         
         \begin{th}[\cite{Nak3}]\label{relative criterion} 
         For $\chi\in Z^1(G, {\rm U}(R))$, the $R^G$-module $R_\chi$ is $R^G$-free
         of rank one if and only if the following conditions are satisfied:
         \begin{itemize}
         \item[(i)] $\dim {\mathcal Q}(R^G)\otimes_{R^G} R_\chi = 1$.
         \item[(ii)] There exists a nonzero element $f \in R_\chi$ satisfying
         \begin{equation} \forall {\frak q} \in {\rm Ht}_1(R^G) ~\Rightarrow~
         \exists {\frak P}\in X_{\frak q}(R) ~\mbox{such that}~
         {\rm v}_{R, {\frak P}}(f) <  {\rm e} ({\frak P}, {\frak q})\label{D5}
          \end{equation}
         \end{itemize}
         If these 
         conditions are satisfied, $R_\chi = R^G \cdot f$
         for any nonzero element $f \in R_\chi$ such that (\ref{D5}) holds.
         
         \end{th}

         \begin{cor}\label{criterion on freeness} 
          Let  $\chi$ be a cocycle in $Z^1(G, {\rm U}(R))_R$. 
        Then
          $R_\chi$ is $R^G$-free if and only if 
          \begin{equation} {\rm div}_R(f)\equiv
           D(\chi) ~~{\rm mod}\bigoplus_{{\frak P}
          \in {\rm Ht}_1^{[2]}(R, R^G)} \mbox{\boldmath $Z$}\cdot
          {\rm div}_R({\frak P})\label{D6}
          \end{equation} 
          for some nonzero $f\in R_\chi$. Especially 
          in the case where  $\chi \in
          -Z^1_R(G, {\rm U}(R))_{(2)}$,  $R_\chi \cong R^G$ as $R^G$-modules 
           if and only if $D(\chi) = {\rm div}_R (f)$ for some nonzero $f\in R_\chi$.  
           Moreover $R_\chi = R^G\cdot f$ holds, 
           in the both  cases where these 
           conditions are
          satisfied.
          \end{cor}
          
          \proof 
          Suppose that the congruence (\ref{D6}) holds 
          for some nonzero $f\in R_\chi$. Then the condition (\ref{D5})
          holds for $f$  and hence, by Theorem \ref{relative criterion} we see that
          $R_\chi \cong R^G$ as $R^G$-modules.  
          Conversely suppose that $R_\chi$ is $R^G$-free. Then by Theorem
          \ref{relative criterion} we can choose
          a nonzero  element from $f\in R_\chi$ in such a way that
          $f$ satisfies the condition (\ref{D5}). By the definition of $D(\chi)$ (cf.
          $D_2(f)$ defined by (\ref{D2})) 
          we must have
          $$D(\chi) = \sum_{{\frak P}\in {\rm Ht}_1^{(0)}(R, R^G)\cup
          {\rm Ht}_1^{(1)}(R, R^G)}
          {\rm v}_{R, {\frak P}}(f) \cdot {\rm div}_R({\frak P}),$$
          which shows that the congruence (\ref{D6}) holds.

          The remainder of the assertion follows easily from Theorem
          \ref{relative criterion} and these observations.
          \qqed
          
          \medskip

     Let ${\rm coh}(\chi)$ denote  the {\it cohomology class} of
     a $1$-cocycle $\chi \in Z^1(G, {\rm U}(R))$.  This  induces a
     homomorphism ${\rm coh} : Z^1(G, {\rm U}(R)) \to H^1(G, {\rm U}(R))$. 
     The zero element of the additive group $H^1(G, {\rm U}(R))$
          is denoted  by ${\rm coh}(\theta)$ (recall that  $\theta$ is
          the trivial $1$-cocycle in $Z^1(G, {\rm U}(R))$). 
     For a
          $G$-invariant principal divisor $D = {\rm div}_R(h)$ for 
          some $h \in {\mathcal Q}(R)$,  define the cohomology class
          ${\rm coh}(D) \in H^1(G, {\rm U}(R))$ which is the
          class of the $1$-cocycle $G \ni \sigma \to \sigma(h)/h \in {\rm U}(R)$. 
          Clearly ${\rm coh}(D)$ does not depend on the choice of $h$. 
        
         For $\chi \in Z^1 (G, {\rm U}(R))$,
             put $\mbox{\boldmath $N$}\cdot\chi:=\{n\chi \in Z^1(G, {\rm U}(R))  \mid
             n \in \mbox{\boldmath $N$}\}$. 
          Recall that  ${\rm tor}(A)$ stands for
           the {\it torsion part} of an abelian group $A$.

         \begin{pr}\label{prop1} 
          Let  $\chi$ be a cocycle in $Z^1_R(G, {\rm U}(R))_{\rm e} 
       \cap\left(-Z^1_R(G, {\rm U}(R))_{(2)}\right)$ such that 
       $\mbox{\boldmath $N$}\cdot\chi \subseteq Z^1(G, {\rm U}(R))_R$. 
             Then the following conditions are equivalent:
             \begin{itemize}
             \item[(i)] $D(\chi)$ is  a principal divisor satisfying the 
             conditions as follows: 
             \begin{itemize}
             \item[(a)]
             ${\rm coh}(\chi) - {\rm coh}(D(\chi)) \not\in {\rm tor}(H^1(G, {\rm U}(R)))
             \backslash \{{\rm coh}(\theta)\}$. 
             \item[(b)]   
             $R_{m\chi} \cong R^G$ as $R^G$-modules  for some $m\in \mbox{\boldmath $N$}$.
              \end{itemize}
             \item[(ii)] For any $n \in \mbox{\boldmath $N$}$,
             $R_{n\chi} \cong R^G$ as $R^G$-modules.
             \item[(iii)]
             $R_\chi \cong R^G$ as $R^G$-modules.
             \end{itemize}

           \end{pr} 
           
           \proof The implication $(ii) \Rightarrow (i)$ follows immediately
           from Corollary \ref{criterion on freeness} (in this case, 
           ${\rm coh}(\chi) = {\rm coh}(D(\chi)$)). 
           
           $(i) \Rightarrow (iii):$ ~Let $f$ be a nonzero element of $R$
           such that ${\rm div}_R(f) = D(\chi)$ and $g$ a nonzero element
           of $R_\chi$. Then, as $R\cdot f \supseteq R\cdot g$ (cf. the inequality
           in {\it (iii)} of Lemma \ref{characteristic divisor}),
           we express $g = f \cdot h$ for some $h\in R$. 
           Since $m\chi \in - Z^1_R(G. {\rm U}(R))_{(2)}$ (cf. Lemma \ref{semigroup}),
           by the choice of $\chi$, {\it (v)} of  Lemma \ref{characteristic divisor}
            and Corollary \ref{criterion on freeness}, we see that
            $$D(m \chi) = {\rm div}_R(f^m) ={\rm div}_R(w)$$ for some
            nonzero element $w\in R_{m  \chi}\cong R^G$. Here note
            $R_{m \chi} = R^G\cdot w$ (cf. the last statement of
            Corollary \ref{criterion on freeness}). 
            So, let $u\in {\rm U}(R)$
            be the unit satisfying $f^m = w \cdot u$. As $R_{m \chi}
             \ni g^m$, we have $$f^m \cdot h^m = g^m = w\cdot v = f^m  
            \cdot u^{-1}\cdot v$$
            for some $v\in R^G$, which requires 
            $${\rm coh}(\chi) - {\rm coh}(D(\chi)) = {\rm coh}({\rm div}_R(h)) \in 
            {\rm tor}(H^1(G, {\rm U}(R))).$$
             Since   $G \ni \sigma \to \sigma (h)\slash h
             \in {\rm U}(R)$ is a coboundary (cf. {\it (a)} of {\it (i)} ), we see     
             $h\cdot t\in R^G$ for a unit
             $t\in {\rm U}(R)$, which implies  $t^{-1} \cdot f \in R_\chi$. 
            Hence the assertion of {\it (iii)} follows from 
            Corollary \ref{criterion on freeness}.

            $(iii) \Rightarrow (ii):$ ~By Corollary \ref{criterion on freeness}, we can choose
            a nonzero element $f\in R_\chi$ in such a way that
            $D(\chi) = {\rm div}_R(f)$. For any $n \in \mbox{\boldmath $N$}$,
            from Lemma \ref{characteristic divisor}
            we infer that  $D(n \chi ) = {\rm div}_R(f^n)$, which
            implies that $R_{n \chi} \cong R^G$ as $R^G$-modules.
            \qqed
            
            \begin{cor}\label{order of characteristic divisor} 
             Under the same circumstances as in Proposition \ref{prop1},  
            suppose that there exists a number $m\in\mbox{\boldmath $N$}$
            such that $R_{m \chi}\cong R^G$ as $R^G$-modules. Then 
            the divisor class $[D(\chi)]$ in ${\rm Cl}(R)$ is a torsion element. 
            Suppose moreover that
            \begin{eqnarray}
            {\rm coh}(i\chi) - {\rm coh}(D(i\chi)) \not\in {\rm tor}(H^1(G, {\rm U}(R)))
             \backslash \{{\rm coh}(\theta)\}\label{D condition}\end{eqnarray}
             for any $i \in \mbox{\boldmath $N$}$ such that
             $D(i\chi)$ is  a principal divisor. 
              Then
           the following equality
            holds; $${\rm ord}( [D(\chi)]) = \min \left\{ q\in \mbox{\boldmath $N$}~
           \left\vert ~R_{q\chi}\cong R^G~\mbox{as $R^G$-modules}\right.\right\}$$
            where ${\rm ord}( [D(\chi)])$ is the order of $[D(\chi)]$
            in the group ${\rm Cl}(R)$. 
            \end{cor}
           
           \proof  Since $m\chi \in - Z^1_R(G. {\rm U}(R))_{(2)}$, 
           by the assumption that $R_{m\chi} \cong R^G$ 
          and Corollary \ref{criterion on freeness}, we see  
          $D(m\chi) = {\rm div}_R(g)$
         for some $g \in R_{m\chi}$, which shows
         $m \cdot [D(\chi)] = 0$ 
        (cf.  {\it (v)} of Lemma
            \ref{characteristic divisor}). We similarly have
            $${\rm ord}( [D(\chi)]) \leqq \min \left\{ q\in \mbox{\boldmath $N$}~
            \left\vert ~ R_{q\chi}\cong R^G~\mbox{as $R^G$-modules}\right.\right\}.$$
        Let $k$ denote the right hand side of this inequality.
        As $R_{k\chi} \cong R^G$, we have an element  $h\in R_{k\chi}$
        satisfying $D(k\chi) = {\rm div}_R(h)$ (cf. Corollary \ref{criterion on freeness}). 
         Let $t\in \mbox{\boldmath $N$}$ be a common multiplier of $k$ and
         ${\rm ord}( [D(\chi)])$. Put $\psi = {\rm ord}([D(\chi)])\chi$.  
         Then,  by {\it (v)} of Lemma
            \ref{characteristic divisor},  $D(t\chi) = (t/k)\cdot D(k\chi)$.  
        Clearly $D(t\chi) = {\rm div}_R(h^{t/k})$
         and $h^{t/k}\in R_{t\chi}$.  
          This implies that $R_{(t/{\rm ord}([D(\chi)]))\psi}
         \cong R^G$ as $R^G$-modules (cf. Corollary \ref{criterion on freeness}). 
         Then from the condition (\ref{D condition}) and the 
          implication $(i) \Rightarrow (iii)$ of Proposition \ref{prop1}
          we infer that $R_{\psi}\cong R^G$ as $R^G$-modules, which completes the proof. 
           \qqed

         \begin{rem}\label{zero cohomology} ~\rm    
         Define $$Z^1(G, {\rm U}(R))_{R, 0} :=
         \{\delta \in Z^1(G, {\rm U}(R))_R \mid D(\delta ) = 0\}.$$
        For $\chi\in Z^1(G, {\rm U}(R))_R$ 
        such that  $D(i\chi)$ is principal,  we have a cocycle 
                $\psi\in Z^1(G, {\rm U}(R))^R$ satisfying 
                $D(i\chi) = {\rm div}_R(g) = D(\psi)$ for some $g \in R_\psi$ and 
                $D(i\chi- \psi) = 0$ (cf. {\it (vii)} of Lemma \ref{characteristic divisor}). 
                Since $${\rm coh}(i\chi) - {\rm coh}(D(i\chi))=
                {\rm coh}(i\chi) - {\rm coh}(D(\psi)) ={\rm coh}(i\chi -\psi),$$ 
          the condition (\ref{D condition})
         holds if 
         the following equality  holds :
         \begin{eqnarray} &\left\{ {\rm coh}(\delta) \mid \delta \in \left(\mbox{\boldmath
         $N$}\cdot\chi +(- Z^1(G, {\rm U}(R))^R)\right) 
            \bigcap Z^1(G, {\rm U}(R))_{R, 0}\right\}\nonumber\\  &
             \cap  {\rm tor}(H^1(G, {\rm U}(R)))
             = \{{\rm coh}(\theta)\}.  \label{torsionfree}
             \end{eqnarray}
               In general we denote by $H^1(G, {\rm U}(R))_{R, 0}$
             the image of $Z^1(G, {\rm U}(R))_{R, 0}$ under the canonical
             map ${\rm coh} : Z^1(G, {\rm U}(R)) \to H^1(G, {\rm U}(R))$. 
             Clearly  $B^1(G, {\rm U}(R)) \subseteq Z^1(G, {\rm U}(R))_{R, 0}$. 
             If the equality
             $$
             H^1(G, {\rm U}(R))_{R, 0} \cap {\rm tor}(H^1(G, {\rm U}(R))) =\{{\rm coh}(\theta)\}
             $$
             holds, then (\ref{torsionfree}) is satisfied. 
         \end{rem}

            \begin{lem}\label{l1}  
           For a cocycle $\chi \in Z^1(G, {\rm U}(R))_R$,
            the $R^G$-module $R_\chi$ is $R^G$-isomorphic to an integral ideal $I$ of $R^G$, which
            is unique up to a multiplication
            of a non-zero
           element  of ${\mathcal Q}(R^G)$.
           Moreover we have an $R^G$-isomorphism
           $d_{(R^G, R)}(R_\chi)\rightarrow
            d_{R^G}(I)$ whose restriction induces $R_\chi \cong I$
          and the divisor class $[d_{R^G}(I)]$
            of $d_{R^G}(I)$
           in ${\rm Cl} (R^G)$ 
           is uniquely determined by $\chi$.  
            \end{lem}
            
            \proof For a nonzero $f \in R_\chi$ we see $R_\chi \cong (1/f)R\cap {\mathcal Q}(R^G)$
            as $R^G$-modules (cf. Lemma 3.1 of \cite{Nak3}) and by Proposition \ref{n1} see that
             the $R^G$-module $(1/f)R\cap {\mathcal Q}(R^G)$
            is a fractional ideal of $R^G$. The assertion  follows immediately from this. \qqed
            
            \begin{df}
            In the circumstances as in Lemma \ref{l1} we denote by $[ R_\chi]$
            the divisor class $[d_{R^G}(I)]\in {\rm Cl}(R^G)$ of $d_{R^G}(I)$ 
            for $\chi \in Z^1(G, {\rm U}(R))_R$, 
            where the ideal $I$ of $R^G$ is defined  in Lemma \ref{l1}. 
            \end{df}

          \begin{pr}\label{prop2} There is a canonical embedding
           $d_{(R^G, R)}(R_\chi) \subseteq R$  
          for any cocycle  $\chi \in Z^1(G, {\rm U}(R))^R\cap Z^1_R(G, {\rm U}(R))_{(2)}$. 
          \end{pr}
          
          \proof Let $\chi$ be a cocycle in $Z^1_R(G, {\rm U}(R))_{(2)}$. Then, 
          for any ${\frak P}\in {\rm Ht}_1(R)$
          such that ${\rm ht} ({\frak P}\cap R^G) \geqq 2$, 
          we can choose a nonzero element $g$ from $R_{-\chi}$ in such a way that
          $g \not\in {\frak P}$.  Since $(R^G)_{{\frak P}\cap R^G}$ is a 
          Krull domain, we see
          \begin{eqnarray*}
          (R^G)_{{\frak P}\cap R^G} &=& \bigcap_{{\frak q}\in{\rm Ht}_1(R^G), {\frak q}\subseteq
          {\frak P}} (R^G)_{{\frak q}} \\ &\supseteq& 
          \bigcap_{{\frak q}\in{\rm Ht}_1(R^G), {\frak q}\subseteq 
          {\frak P}}  (R^G)_{{\frak q}}\otimes_{R^G} (g\cdot R_\chi) 
                 \supseteq  
             g\cdot d_{(R^G, R)}({R_\chi}).
          \end{eqnarray*}   Thus by the choice of $g$, one sees
          $ g\cdot R_{\frak P} = R_{\frak P} \supseteq
          (R^G)_{{\frak P}\cap R^G} \supseteq g\cdot d_{(R^G, R)}({R_\chi})$, which requires
          $R_{\frak P} \supseteq d_{(R^G, R)}({R_\chi})$ for any $\frak P\in {\rm Ht}_1(R)$
          such that ${\rm ht} ({\frak P}\cap R^G) \geqq 2$. 
          \par In general, we have
          \begin{eqnarray}  
          R_{\frak P} &\supseteq&  R_{\frak P}\otimes_R (R\cdot R_\chi) \nonumber\\
                                &\supseteq& (R\backslash {\frak P})^{-1} R_\chi
                                \left(= \left
                                \{\left. \frac{a}{b}~
                                \right
                                \vert  ~a\in R_\chi, ~b\in R\backslash {\frak P}\right\}\right)
                                \nonumber\\
                                  &\supseteq& R^G_{{\frak P}\cap R^G}\otimes_{R^G} R_\chi
                                \label{eqn:1}
          \end{eqnarray}
          for a prime ideal  ${\frak P}$ of $R$. 
          Here we identify $R^G_{{\frak P}\cap R^G}\otimes_{R^G} R_\chi$ (resp. 
          $R_{\frak P}\otimes_R (R\cdot R_\chi)$) with 
          $(R^G\backslash ({\frak P}\cap R^G))^{-1} R_\chi$ 
           (resp. $(R\backslash {\frak P})^{-1} R\cdot R_\chi$) in ${\mathcal Q(R)}$.
            Since $X_{\frak q}(R) \not= \emptyset$ for all ${\frak q}\in {\rm Ht}_1(R^G)$,  
        by (\ref{eqn:1}) we see that
         \begin{eqnarray}
              &\bigcap_{{\frak q}\in {\rm Ht}_1(R^G)} 
              \left( (R^G)_{{\frak q}}\otimes_{R^G} R_\chi\right)
              = \bigcap_{{\frak P}\in {\rm Ht}_1^{(1)}(R, R^G)} 
              \left((R^G)_{{\frak P}\cap R^G}\otimes_{R^G} 
              R_\chi\right)\nonumber \\
              &= \bigcap_{{\frak P}\in {\rm Ht}_1^{(1)}(R, R^G)} 
              \left( (R^G)_{{\frak P}\cap R^G}\otimes_{R^G} 
              R_\chi\right)\cap \bigcap_{{\frak P}\in {\rm Ht}_1^{(0)}(R, R^G)}
              \left((R^G)_{{\frak P}\cap R^G} \otimes_{R^G} R_\chi\right)\nonumber\\
              &\subseteq
               \bigcap_{{\frak P}\in {\rm Ht}_1(R)
               \setminus{\rm Ht}_1^{[2]}(R, R^G)} \left( R_{{\frak P}}\otimes_{R} 
              R \right)
        \nonumber.
          \end{eqnarray}
          These localizations are regarded as modules of fractions in ${\mathcal Q}(R)$
          and intersections are defined in ${\mathcal Q}(R)$.
         Thus
          we must have
         \begin{eqnarray}
         d_{(R^G, R)}({R_\chi})=
         \bigcap_{{\frak q}\in {\rm Ht}_1(R^G)} \left( R^G_{{\frak q}}\otimes_{R^G} R_\chi\right)
         \subseteq  \bigcap_{{\frak P}\in {\rm Ht}_1(R)} \left( R_{{\frak P}}\otimes_R R\right)
         = R \nonumber 
         \end{eqnarray}
          because $R$ is a Krull domain. 
          \qqed
          
          \begin{cor}\label{main corollary}
          Let  $\chi$ be a cocycle in $Z^1_R(G, {\rm U}(R))_{\rm e} 
       \cap Z^1_R(G, {\rm U}(R))_{(2)}$.  
         Then, for a natural number   $n$,
        $R_{n\chi} \cong R^G$ as $R^G$-modules
        if and only if ~$n \cdot [R_\chi] = 0$ in ${\rm Cl} (R^G)$. 
         \end{cor}
          
          \proof  By the choice of $\chi$, we   see
          $$\mbox{\boldmath $N$} \cdot \chi \subseteq Z^1(G, {\rm U}(R))^R
          \cap \bigl( -Z^1(G, {\rm U}(R))^R\bigr) \subseteq  Z^1(G, {\rm U}(R))_R.$$
          Let $f$ be a nonzero element of $R_\chi$.  Applying  
          Corollary \ref{power of ideal}
          to the ideal $\displaystyle R\cdot({1}/{f})$ with the aid of {\it (ii)}  
          of Lemma \ref{characteristic divisor},
          we see \begin{equation}
          \displaystyle {\rm div}_{R^G}\Bigl( R\cdot\frac{1}{f^n}\cap {\mathcal Q}(R^G)
          \Bigr) = \displaystyle n\cdot 
          {\rm div}_{R^G}\Bigl( R\cdot\frac{1}{f}\cap {\mathcal Q}(R^G)
          \Bigr)\label{eqn8}\end{equation}
          for any $n \in \mbox{\boldmath $N$}$ (cf. Lemma 3.1 of \cite{Nak3}). 
       From 
          Lemma \ref{l1},  the 
          identity (\ref{eqn8}) can be deduced to  $[R_{n \chi}] = n\cdot [R_{\chi}]$
          in ${\rm Cl}(R^G)$. By Lemma \ref{semigroup}, 
           Lemma \ref{divisorialization} and Proposition 
          \ref{prop2}, we immediately see that $R_{n \chi} \cong R^G$ if and only if
          $[R_{n\chi}] = 0$ in ${\rm Cl}(R^G)$, which implies the assertion. \qqed
          
          \medskip

              Combining Corollary \ref{order of
          characteristic divisor} with Corollary \ref{main corollary}, we immediately have
          
         \begin{th}\label{main Krull case}
          Let  $\chi$ be a cocycle in $Z^1_R(G, {\rm U}(R))_{\rm e} 
       \cap\widetilde{ Z}^1_R(G, {\rm U}(R))$.  
        Suppose that the equality
       (\ref{D condition}) holds 
       for any $i \in \mbox{\boldmath $N$}$ such that
             $D(i\chi)$ is  a principal divisor. 
               If $[R_\chi]\in {\rm tor} ({\rm Cl}(R^G))$,
       then $${\rm ord} ([R_\chi]) ~\mbox{in ${\rm Cl}(R^G)$}=
        {\rm ord} ([D(\chi)]) ~\mbox{in ${\rm Cl}(R)$} ~(< \infty)$$
        which  is equal to \(\min \left\{ q\in \mbox{\boldmath $N$}
            \mid R_{q \chi}\cong R^G~\mbox{as $R^G$-modules}\right\}\).
            \qqed
           \end{th}

          \begin{df}\label{theta}
          \rm  Let ${\rm UrCl}(R, G)$ denote the
          subgroup of ${\rm Cl}(R)$ generated by 
          $$\Bigl\{[D(\chi)] ~\Big\vert~\chi\in Z^1_R(G, {\rm U}(R))_{\rm e} 
       \cap\widetilde{Z}^1_R(G, {\rm U}(R))\Bigr\}$$ where
          $[D(\chi)]$ denotes the divisor class of $D(\chi) \in {\rm Div} (R)$. 
         Define 
           $\widetilde{{\rm Cl}}(R, G)$  to be the  subgroup $$\left< \Bigr\{[R_\chi] 
           ~\Bigm\vert~
            \chi \in Z^1_R(G, {\rm U}(R))_{\rm e} 
       \cap\widetilde{Z}^1_R(G, {\rm U}(R))\Bigr\}\right>$$
       of ${\rm Cl} (R^G)$. 
            \end{df}
        
        Recall that  ${\rm exp} (N)$ denotes  the {\it exponent} of a group $N$. 
          The next result follows 
         from Remark \ref{zero cohomology} and Theorem \ref{main Krull case}.

           \begin{pr}\label{final}
         Suppose that the subset 
            $${\rm coh}\Bigl(Z^1(G, {\rm U}(R))_{R,{0}}\cap
            Z^1_R(G, {\rm U}(R))_{\rm e}\cap
           \widetilde{Z}^1_R(G, {\rm U}(R))\Bigr)$$ 
           of  $H^1(G, {\rm U}(R))$
           does not contain a non-trivial torsion element of $H^1(G, {\rm U}(R))$.  
           Suppose that 
             $\widetilde{\rm Cl}(R, G)$ is a 
            torsion group. If one of  ${\rm exp} ({\rm UrCl}(R, G))$ 
            and ${\rm exp} (\widetilde{{\rm Cl}}(R, G))$ is finite, then we have 
            ${\rm exp} ({\rm UrCl}(R, G)) = {\rm exp} (\widetilde{{\rm Cl}}(R, G))$
      \end{pr}
      
      \proof Let $\chi$ be a cocycle in $Z^1_R(G, {\rm U}(R))_{\rm e}\cap
           \widetilde{Z}^1_R(G, {\rm U}(R))$ and 
           $i$ a natural number such that $D(i\chi)= {\rm div}_R(g)$ for a 
           nonzero $g \in R$. Let $\psi$ be a cocycle of $G$ in ${\rm U}(R)$ 
           defined by  $R_\psi \ni g$. Then $\psi \in Z^1_R(G, {\rm U}(R))_R$,
           $D(\psi) = D(i\chi)$, $R_\psi = R^G \cdot g$,
           $i\chi - \psi \in Z^1_R(G, {\rm U}(R))_R$ and $D(i\chi -\psi) = 0$
           (cf. {\it (vii)} of Lemma \ref{characteristic divisor}).  
      Since 
           $R_{i\chi} = R_{i\chi - \psi}\cdot g$ (cf. {\it (iv)} of Lemma
           \ref{characteristic divisor}), $R_{-i\chi}\cdot g
           \subseteq R_{-i\chi +\psi}$
                and $i\chi \in 
           Z^1_R(G, {\rm U}(R))_{\rm e}\cap \widetilde{Z}^1_R(G, {\rm U}(R))$, 
           by {\it (v)} of Lemma \ref{characteristic divisor} and  (\ref{D2}) we have  
            $$i\chi -\psi\in Z^1_R(G, {\rm U}(R))_{\rm e}\cap\widetilde{Z}^1_R(G, {\rm U}(R)).$$ 
            Thus the condition (\ref{D condition}) holds for $i\chi$ and 
           Theorem \ref{main Krull case} can be applied to $\chi$.  
            By our assumption, 
       ${\rm UrCl}(R, G)$ and $\widetilde{{\rm Cl}}(R, G)$ are 
      torsion groups. The exponent of a torsion abelian group
      generated by a subset $\Xi$ is determined by 
      the ideal $\displaystyle \cap _{a\in \Xi} ~{\rm ord}(a) \mbox{\boldmath $Z$}$
      of $\mbox{\boldmath $Z$}$.
   The assertion follows from this observation and Theorem \ref{main Krull case}. \qqed
      
            \begin{th}\label{finiteness of class group}  Suppose that 
                     \begin{equation}
            {\rm coh}\bigl(Z^1_R(G, {\rm U}(R))_{\rm e}\cap
           \widetilde{Z}^1_R(G, {\rm U}(R)\bigr)\subseteq
            \sum_{\lambda \in \Lambda}\mbox{\boldmath $Z$}_0
            \cdot {\rm coh}(\lambda)\label{semigroup condition}\end{equation}
            for a non-empty finite subset $\Lambda$ of $Z^1(G, {\rm U}(R))^R$. 
            Suppose that   
             $\widetilde{\rm Cl}(R, G)$ is a 
            torsion group. Then it is a finite group.  Moreover if
             $${\rm coh}\Bigl(Z^1(G, {\rm U}(R))_{R, {0}}\cap
            Z^1_R(G, {\rm U}(R))_{\rm e}\cap
           \widetilde{Z}^1_R(G, {\rm U}(R))\Bigr)$$ in $H^1(G, {\rm U}(R))$
           does not contain a non-trivial torsion element, we have
            $${\rm exp} ({\rm UrCl}(R, G)) = {\rm exp} (\widetilde{{\rm Cl}}(R, G)) <\infty.$$
            \end{th}
     
     \proof  Let $\Lambda = \{\lambda_1, \dots, \lambda_n\}$ and choose
     a nonzero $g_i$ from $R_{\lambda_i}$ for any $1 \leqq i \leqq n$. By (\ref{semigroup
     condition}), for $\chi \in Z^1_R(G, {\rm U}(R))_{\rm e}\cap
           \widetilde{Z}^1_R(G, {\rm U}(R))$, 
            we can choose $j_i\in \mbox{\boldmath $Z$}_0$ and $u \in
           {\rm U}(R)$ in such a way that
           $\chi - \sum_{i=1}^n j_i \lambda_i$ is a $1$-coboundary
           of $G$ defined by $u^{-1}$. Then $$d_{(R^G, R)}(R_\chi)
           \cong d_{(R^G, R)}\Bigl(\frac{u}{~\prod g_i^{j_i}~}R\cap {\mathcal Q}(R^G)\Bigr)
           = d_{(R^G, R)}\Bigl(\frac{1}{~\prod g_i^{j_i}~}R\cap {\mathcal Q}(R^G)\Bigr)$$
           as $R^G$-modules (cf. Lemma 3.1 of \cite{Nak3}). 
           Consequently  $\widetilde{{\rm Cl}}(R, G)$
           is a subgroup of the image of $D_{(R, R^G)}(g_1, \dots, g_n)$
           (for definition,  cf. Lemma \ref{finite lemma}) 
           under the canonical homomorphism ${\rm Div}(R^G) \to {\rm Cl}(R^G)$.
           By Lemma \ref{finite lemma} the group $\widetilde{{\rm Cl}}(R, G)$ is 
           finitely generated.  The remainder of the
            assertion of this theorem follows  
           from Proposition \ref{final}. \qqed
       
         \small \section{Affine normal varieties with group actions} \normalsize

        Hereafter let $G$ be an affine algebraic group over an algebraically closed
        field $K$ of characteristic $p \geqq 0$. 
        Let $(X, G)$ be a regular  
        action of $G$ 
        on an {\it affine normal 
         variety} $X$ defined over  $K$.
         Let ${\frak X}(G)$ be the additive group of all rational characters of
          $G$ and for a rational $G$-module $W$ and $\chi\in {\frak X}(G)$, denote by
          $W_\chi$ the subspace $\{ x \in W \mid \sigma(x) = 
          \chi (\sigma) x ~(\forall \sigma \in G)\}$
          of $W$ consisting of all vectors of weight $\chi$.  Put
          ${\frak X}(G)^W :=\{\chi \in {\frak X}(G) \mid W_\chi \not=\{0\}\}$.  
          The cocycles
          of $G$ on ${\rm U}({\mathcal O}(X))$ defined by 
          regular functions  are identified with  rational characters:
          If  $G$ is connected, by a result of 
          M. Rosenlicht we see  $Z^1(G, {\rm U}({\mathcal O}(X)))^{{\mathcal O}
          (X)} = {\frak X}(G)^{{\mathcal O}(X)}$ (e.g.,  \cite{Bass2, Magid}), since $X$
          is normal.

         For any ${\frak P}\in {\rm Ht}_1^{(1)}
        ({\mathcal O}(X), {\mathcal O}(X)^G)$, put $${\mathcal I}_G({\frak P}): =
        \{\sigma \in G \mid \sigma (x)~ \equiv ~x ~{\rm mod} ~{\frak P} ~~(\forall x
        \in {\mathcal O}(X)) \}$$
        which is called the {\it inertia group at} ${\frak P}$ under the action 
        of $G$.  Clearly
        ${\mathcal I}_G({\frak P})$ contains the ineffective
        kernel $  {\rm Ker} (G \to {\rm Aut} (X))$ of the action
         $(X, G)$ which is denoted by 
        $L_{(X, G)}$. 
        \begin{df}\label{def of pseudo-reflection group}~ \rm The {\it
             pseudo-reflection group of the action} $(X, G)$, denoted by  
           ${\frak  R}(X, G)$,  is defined to be the subgroup of $G$
           generated  by $
             {\mathcal I}_G({\frak P})$'s for all
             ${\frak P}\in {\rm Ht}_1^{(1)}({\mathcal O}(X), 
           {\mathcal O}(X)^G)$.  
         Let $\widetilde{{\frak  R}}(X, G)$  be the subgroup of $G$
          generated by $L_{(X, G)}$ and   ${\mathcal I}_G({\frak P})$'s
          for all \({\frak P}\in {\rm Ht}_1^{(1)}({\mathcal O}(X), 
           {\mathcal O}(X)^G)
           \) such that ${\frak P}$ are non-principal.
          This is called the {\it  non-principal  pseudo-reflection subgroup
          of the action} $(X, G)$. 
          \end{df} 
          
          Obviously $G \rhd {\frak  R}(X, G) \rhd \widetilde{{\frak  R}}(X, G) \rhd
          L_{(X, G)}$. 
          By Sect. 1 of \cite{Nak3} the group ${\frak  R}(X, G)\vert _{X}$ \linebreak 
          ($\cong {\frak  R}(X, G)/ L_{(X, G)}$) 
             of  restrictions  to $X$
             is finite under the assumption that $G^0$ is linearly reductive. 
           In the
           case where ${\mathcal O}(X)^G$ is finitely generated  over $K$, 
           let $X\dslash G$ denote the affine variety defined by ${\mathcal O}(X)^G$
           and $\pi _{X, G} : X \to X\dslash G$ denote
           the quotient morphism defined by 
           the inclusion ${\mathcal O}(X)^G \to {\mathcal O}(X)$.
           
            \begin{df} \rm
           The action $(X, G)$ admitting its quotient $X\dslash G$, i.e., 
           ${\mathcal O}(X)^G$ is finitely generated over $K$,  
           is defined  to be {\it cofree} (resp. {\it equidimensional}\/),
           if ${\mathcal O}(X)$ is a free ${\mathcal O}(X)^G$-module
           (resp. $\pi _{X, G} : X \to X\dslash G$ is equidimensional).
           Especially in the case where 
            $G$ is linearly reductive, $(X, G)$ is defined to be
           {\it isobaric cofree} if the  $\psi$-isotypical 
           component ${\mathcal O}(X)_{[\psi]}
           := \sum_{V\subseteq {\mathcal O}(X), V\cong \psi}
           V\subseteq {\mathcal O}(X)$ of the rational $G$-module
           ${\mathcal O}(X)$ is zero or ${\mathcal O}(X)^G$-free,
           for any irreducible representation $\psi$
           of $G$. 
            \end{df}
            
             For  a conical action $(X, G)$ of a conical $X$ with
           a linearly reductive $G$, {\it $(X, G)$ is cofree if and only if
           it is isobaric cofree.} Here $X$ is said to be {\it conical}, 
           if the coordinate ring ${\mathcal O}(X)$ is a $\mbox{\boldmath $Z$}_0$-graded
           algebra defined over ${\mathcal O}(X)_0 = K$ and  $(X, G)$ is said to be
           {\it conical}, moreover if the action of 
           $G$ preserves each homogeneous part of ${\mathcal O}(X)$. 
          We  immediately  have
           
             \begin{lem}\label{no blowing up}
              Suppose that ${\mathcal O}(X)^G$ is finitely generated over $K$. 
            If ${\mathcal O}(X)^G \to {\mathcal O}(X)$ is no-blowing-up of codimension one
            (e.g., PDE in p. 30 of \cite{Fossum}) (especially if $(X, G)$ is equidimensional), then
            $Z^1(G, {\rm U}(R))^R \cap (-Z^1(G, {\rm U}(R))^R ) =
            \widetilde{Z}^1_R(G, {\rm U}(R)). ~ \Box $
            
            \end{lem}

          \begin{pr}\label{generalized reflection}
          Suppose that $G^0$ is linearly 
          reductive. Let $N$ be a normal closed subgroup of $G$
          such that $N\vert_X$ is finite. Then:
          \begin{itemize}
          \item[{\it (i)}]
           ${\frak R}(X\dslash N, G) = N\cdot {\frak R}(X, G)$ and 
           $N\lhd\widetilde{{\frak R}}(X\dslash N, G)\lhd N\cdot {\frak R}(X, G)$. 
           \item[{\it (ii)}]
           If both $X$ and $(X, G)$ are conical and the order 
           $\vert (N\cdot {\frak R}(X, G))/\widetilde{{\frak R}}(X\dslash N, G)\vert$
           is a unit in $K$, then the natural  action 
           $(X\dslash \widetilde{{\frak R}}(X\dslash N, G), N\cdot {\frak R}(X, G))$
           is cofree. 
          
          \end{itemize}  
          
          \end{pr}
          
          \proof {\it (i):}~ Since $N$ is normal in $G$, 
          we immediately see ${\frak R}(X\dslash N, G) \supseteq {\frak R}(X, G).$
          For any  ${\frak p}\in {\rm Ht}_1^{(1)}({\cal O}(X)^N, {\cal O}(X)^G)$ let 
           ${\frak P}$ be a prime ideal of  ${\cal O}(X)$ lying over ${\frak p}$.  
          Put $H = {\mathcal I}_G({\frak p})$. Clearly $ H\supseteq N$, ${\mathcal I}
          _H({\frak P}) = {\mathcal I}_G({\frak P})$ and 
          $H\vert _X$ is finite (cf. Sect. 1 of \cite{Nak2}). From Expos\'e V of 
          \cite{Gr} we infer that the canonical morphisms
          $${\mathcal O}(X)^{N\cdot {\mathcal I}_H({\frak P})} \to 
          {\mathcal O}(X)^{{\mathcal I}_H({\frak P})}$$
           $${\mathcal O}(X)^{H} \to 
          {\mathcal O}(X)^{{\mathcal I}_H({\frak P})}$$
          are \'etale at ${\frak P}\cap {\mathcal O}(X)^{{\mathcal I}_H({\frak P})}$. 
          Consequently the monomorphism $${\mathcal O}(X)^{H} \to 
          {\mathcal O}(X)^{N\cdot {\mathcal I}_H({\frak P})}$$ is unramified 
          at ${\frak p}\cap {\mathcal O}(X)^{N\cdot {\mathcal I}_H({\frak P})}$. 
         If $\sigma\in H$, we see $\sigma \cdot {\mathcal I}_H({\frak P}) \cdot \sigma^{-1}
         =  {\mathcal I}_H(\sigma ({\frak P})) = {\mathcal I}_H(\tau ({\frak P}))
         = \tau \cdot {\mathcal I}_H({\frak P}) \cdot\tau^{-1}$
         for some $\tau \in N$. Thus $N\cdot {\mathcal I}_H({\frak P})$
         is a normal subgroup of $H$ stabilizing 
         ${\frak p}\cap {\mathcal O}(X)^{N\cdot {\mathcal I}_H({\frak P})}$. 
         On the other hand, by the definition of $H$ and Sect. 41 of \cite{LR}, 
         the residue class field of $({\mathcal O}(X)^{N\cdot {\mathcal I}_H({\frak P})})_
          {{\frak p}\cap {\mathcal O}(X)^{N\cdot {\mathcal I}_H({\frak P})}}$
          is purely inseparable over that  of 
          $({\mathcal O}(X)^{H})_{{\frak p}\cap {\mathcal O}(X)^{H}}$, which
          requires that these fields coincide. 
         Applying
         Nakayama Lemma to the finite unramified local morphism
         $$({\mathcal O}(X)^{H})_{{\frak p}\cap {\mathcal O}(X)^{H}}  \to 
          ({\mathcal O}(X)^{N\cdot {\mathcal I}_H({\frak P})})_
          {{\frak p}\cap {\mathcal O}(X)^{N\cdot {\mathcal I}_H({\frak P})}},$$
          we must have $H = N\cdot {\mathcal I}_H({\frak P}) = N\cdot {\mathcal I}_G({\frak P})$,
          which shows the assertion. 
          
          {\it (ii):}~ Let ${\frak p}\in {\rm Ht}_1^{(1)}({\cal O}(X)^N, {\cal O}(X)^G)$
          be a principal ideal such that ${\mathcal I}_G({\frak p})\vert _{X\dslash N}$
          is non-trivial.  Then there exists a {\it homogeneous} element $f_1$ generating 
          principally  ${\frak p}$. Let $\{Kf_1, Kf_2, \dots, Kf_s\}$
          be the $N\cdot {\frak R}(X, G)$-orbit of $Kf_1$ consisting of $s$ $K$-subspaces
        of ${\mathcal O}(X)^N$.  Let ${\frak q}$ be
          a non-principal ideal in 
          ${\rm Ht}_1^{(1)}({\cal O}(X)^N, {\cal O}(X)^G)$.  For any $\tau \in 
          {\mathcal I}_G({\frak q})$,  we easily see that $\tau (\prod_{i=1}^s f_i)
          = \prod_{i=1}^s f_i$. Thus ${\frak p}\cap {\mathcal O}(X)^{\widetilde{\frak R}
          (X\dslash N, G)} = {\mathcal O}(X)^{\widetilde{\frak R}
          (X\dslash N, G)} \prod_{i=1}^s f_i$ and 
          $$(\sigma -1)({\mathcal O}(X)^{\widetilde{\frak R}(X\dslash N, G)})
          \subseteq {\mathcal O}(X)^{\widetilde{\frak R}
          (X\dslash N, G)} \prod_{i=1}^s f_i$$ for any $\sigma \in {\mathcal I}_G({\frak p})$.
          Since $N\cdot {\frak R}(X, G)/ \widetilde{\frak R}(X\dslash N, G)$
          is generated by generalized reflections in ${\rm Aut}
          ({\mathcal O}(X)^{\widetilde{\frak R}
          (X\dslash N, G)})$ in the sense of M. Hochster and J. A. Eagon \cite{HE}, the 
          ${\mathcal O}(X)^{N\cdot {\frak R}(X, G)}$-module
          ${\mathcal O}(X)^{\widetilde{{\frak R}}(X\dslash N, G)}$ is free
          (e.g., Chapitre 5 of \cite{Bourbaki}).  \qqed

          \begin{pr}\label{pseudo-reflection group}
            Suppose that ${\rm char} (K) = p = 0$,  $G^0$ is an algebraic torus, 
          the action $(X, G)$ is stable and $G$ equals to the centralizer $Z_G(G^0)$ of $G^0$
          in $G$.
          Then 
         $${\frak X}(G)\cap
          Z^1_{{\mathcal O}(X)}(G, {\rm U}({\mathcal O}(X)))_{\rm e} = 
          \bigl\{ \chi \in {\frak X}(G) 
          ~\bigm\vert ~ {\mathcal O}(X)_\chi \not= \{0\}, 
          \chi ({\frak  R}(X, G)) = \{1\}  \bigr\}.$$
        
          \end{pr}
          
          \proof  Let  ${\frak P}\in {\rm Ht}_1^{(1)}
        ({\mathcal O}(X), {\mathcal O}(X)^G)$ and $\chi \in  {\frak X}(G)\cap
        Z^1(G, {\rm U}({\mathcal O}(X)))^{{\mathcal O}(X)}$. Suppose that there exists 
        a nonzero
        $f_{\frak P} \in R_\chi$ such that ${\rm v}_{{\mathcal O}(X), {\frak P}}(f_{\frak P})$
        is divisible by the ramification index
         ${\rm e}({\frak P}, {\frak P}\cap {\mathcal O}(X)^G)$. Then, as
         ${\mathcal O}(X)^G_{{\frak P}\cap {\mathcal O}(X)^G}$ is a discrete
         valuation ring, $f_{\frak P} =
          g\cdot  w$
         for some $g\in {\rm U}({\mathcal O}(X)_{{\frak P}})$ and 
           $w\in {\mathcal O}(X)^G_{{\frak P}\cap {\mathcal O}(X)^G}$. For any
           $\tau \in {\mathcal I}_G({\frak P})$, we have $\chi(\tau) g\cdot w 
           = \tau (g) \cdot w$ and hence 
           $\chi (\tau )~ \equiv~ 1~{\rm mod}~{\frak P}$, which induces
           $\chi (\tau ) = 1$. 
           
           Conversely suppose that $\chi ({\mathcal I}_G({\frak P}))= \{1\}$. 
           For any nonzero $f\in R_\chi$,  we see
           \begin{eqnarray}
           {\rm v}_{{\mathcal O}(X), {\frak P}}(f) &=& {\rm e}({\frak P},
           {\frak P}\cap {\mathcal O}(X)^{{\mathcal I}_G({\frak P})})\cdot 
           {\rm v}_{{\mathcal O}(X)^{{\mathcal I}_G({\frak P})}, {\frak P}\cap
           {\mathcal O}(X)^{{\mathcal I}_G({\frak P})}}(f)\nonumber\\ &\equiv& 0~ 
           {\rm mod} ~{\rm e}({\frak P}, {\frak P}\cap {\mathcal O}(X)^G) \nonumber
           \end{eqnarray}
           as ${\rm e}({\frak P},
           {\frak P}\cap {\mathcal O}(X)^{{\mathcal I}_G({\frak P})})
           = {\rm e}({\frak P}, {\frak P}\cap {\mathcal O}(X)^G)$ (cf. \cite{Nak2}).
           This shows the assertion. 
           \qqed
           
           \medskip

            For any $m\in \mbox{{\boldmath
           $N$}}$, 
           put  $${\rm tor}(m, G, H):=
           \{\sigma\in G \mid \sigma ^m \in H\}$$ for a closed subgroup $H$ of $G$.  
           When $G \rhd H$ and $G/H$ is abelian,  ${\rm tor}(m, G, H)$ is a 
           normal closed subgroup of $G$ containing $H$. 
           Suppose that $G$
           is connected. Put
           $${\frak X}(G)_{{\mathcal O}(X)} : =
          \{ \chi \in  {\frak X}(G) \mid {\mathcal O}(X)_\chi \cdot 
          {\mathcal O}(X)_{-\chi} \not= \{0\}\}$$ and define
          ${\mathcal K}(X, G) := \cap _{\chi \in {\frak X}(G)_{{\mathcal O}(X)}} {\rm Ker} (\chi)$. 
          Then  $${\mathcal O}(X)^{{\mathcal K}(X, G)} = K[ 
          {\mathcal O}(X)_\chi \mid \chi \in {\frak X}(G)_{{\mathcal O}(X)}]
          = \bigoplus _{\chi \in {\frak X}(G)_{{\mathcal O}(X)}} {\mathcal O}(X)_\chi$$ and if this
          $K$-algebra is  finitely generated, the induced action 
           $(X\dslash {\mathcal K}(X, G), G)$
          is  stable (cf. \cite{Nak1}). 
          In  the case where  $G$ is an algebraic torus,  
          $(X, G)$ is stable if and only if $X = X\dslash {\mathcal K}(X, G)$.

          Let $G \times {\frak X}(G) \to K^\ast$ denote the canonical pairing,
          where $K^\ast$ denotes  ${\rm U}(K)$.  
          The {\it orthogonal set operation} $\perp_G$ is  defined naturally by this pairing, i.e,
          $Y^{\perp_G} := \bigcap_{\psi\in Y}{\rm Ker}(\psi)$ for a subset
          $Y \subseteq {\frak X}(G)$
          and $N^{\perp_G}:= \left\{\psi \in {\frak X}(G) \mid \psi (N) =\{1\}\right\}$
          for a subset $N\subseteq G$.

          \begin{lem}\label{pairing}
          Suppose that ${\rm char}(K) = p = 0$ and $G$ is a  connected algebraic group.
          Let $\Gamma$ be a closed normal subgroup of $G$ containing ${\mathcal K}(X, G)$. 
           Then:
          \begin{itemize}
          \item[{\it (i)}]  $\Gamma^{\perp_G} = {\frak X}(G)_{{\mathcal O}(X)^\Gamma}$ 
          and $\left({\frak X}(G)_{{\mathcal O}(X)^\Gamma}\right)^{\perp_G} = \Gamma$. 
          \item[{\it (ii)}] For any $m \in \mbox{\boldmath $N$}$, 
          ${\rm tor}(m, G, \Gamma)^{\perp_G} = 
          m \cdot {\frak X}(G)_{{\mathcal O}(X)^\Gamma}
          \subseteq {\frak X}(G)$. 
        
          \end{itemize} 
           Especially, putting $\Omega_{(X, G)} := \widetilde{Z}^1_{{\mathcal O}(X) } 
          (G, {\rm U}({\mathcal O}(X)))^{\perp_G}$,
          we have 
          $$\widetilde{Z}^1_{{\mathcal O}(X) } 
          (G, {\rm U}({\mathcal O}(X)))=
            \widetilde{Z}^1_{{\mathcal O}(X)^{\Omega_{(X, G)}} }
            (G, {\rm U}({\mathcal O}(X)^{\Omega_{(X, G)}})) = 
            {\frak X}(G)_{{\mathcal O}(X)^{\Omega_{(X, G)}}}.$$
           
          \end{lem}
          
          \proof Since the canonical pairing
          $$G/(\cap_{\psi\in {\frak X}(G)}{\rm Ker}(\psi)) \times
          {\frak X}(G/(\cap_{\psi\in {\frak X}(G)}{\rm Ker}(\psi))) \longrightarrow K^\ast$$
          of a diagonalizable group
          is non-degenerate, 
          $\left(\Gamma^{\perp_G}\right)^{\perp_G} = \Gamma$ and 
          $\left(({\frak X}(G)_{{\mathcal O}(X)^\Gamma})^{\perp_G}\right)^{\perp_G}
          = {\frak X}(G)_{{\mathcal O}(X)^\Gamma}$. 
          Especially $${\mathcal K}(X, G)^{\perp_G} = 
          \left(({\frak X}(G)_{{\mathcal O}(X)})^{\perp_G}
          \right)^{\perp_G} = {\frak X}(G)_{{\mathcal O}(X)} 
          \supseteq \Gamma^{\perp_G}$$ and hence for $\psi \in \Gamma^{\perp_G}$
          we have $\{0\} \not=
          {\mathcal O}(X)_{\pm\psi} \subseteq {\mathcal O}(X)^{\Gamma}$,
          which shows $\Gamma^{\perp_G}\subseteq {\frak X}(G)
          _{{\mathcal O}(X)^{\Gamma}}$. The converse inclusion is
          obvious. Thus the assertion in {\it (i)} follows.
          
          {\it (ii) :}  Consider the canonical non-degenerate pairing
          $G/\Gamma \times {\frak X}(G/\Gamma) \to K^\ast$. Then by 
          the definition of ${\rm tor}(m, G, \Gamma)$,
          we have $$\bigl({\rm tor}(m, G, \Gamma)/\Gamma \bigr)^{\perp_{G/\Gamma}}
          = {\rm tor}(m, G/\Gamma, \{\Gamma\})^{\perp_{G/\Gamma}}
          = m\cdot {\frak X}(G/\Gamma)$$ (e.g.,  Chap. 3 of \cite{Sp}). 
          This implies ${\rm tor}(m, G, \Gamma)^{\perp_G} =
          m \cdot \Gamma^{\perp_G}$ and the assertion follows from 
          {\it (i)}. 
          
           Let $\chi\in \widetilde{Z}^1_{{\mathcal O}(X) } 
          (G, {\rm U}({\mathcal O}(X)))$ be a cocycle. For  ${\frak Q}\in
          {\rm Ht}^{[2]}_1({\mathcal O}(X)^{\Omega_{(X, G)}}, {\mathcal O}(X)^G)$, we
          have $\frak P\in X_{\frak Q}({\mathcal O}(X))$ and,  by the choice
          of $\chi$, $${\mathcal O}(X)_{\pm\chi}
          = {\mathcal O}(X)_{\pm\chi}\cap {\mathcal O}(X)^{\Omega_{(X, G)}}
           \not\subseteq {\frak P} \cap {\mathcal O}(X)^{\Omega_{(X, G)}} ={\frak Q},$$ 
           which implies that  $$\widetilde{Z}^1_{{\mathcal O}(X) } 
          (G, {\rm U}({\mathcal O}(X)))\subseteq
            \widetilde{Z}^1_{{\mathcal O}(X)^{\Omega_{(X, G)}} }
            (G, {\rm U}({\mathcal O}(X)^{\Omega_{(X, G)}}))
            \subseteq {\frak X}(G)_{{\mathcal O}(X)^{\Omega_{(X, G)}})}.$$ 
            Here
            $Z^1(G, {\rm U}({\mathcal O}(X)^{\Omega_{(X, G)}}))$
            is naturally regarded as a subgroup of 
            $Z^1(G, {\rm U}({\mathcal O}(X))$. The last assertion can be
            shown by applying  
            {\it (i)} to $\Gamma = \Omega_{(X, G)}$.   
            \qqed


       \small   \section{Equidimensional toric actions} \normalsize

       Denote by   $\Delta_{(X, G)}$ the group  $\left(
         \widetilde{Z}^1_{{\mathcal O}(X)} 
          (G, {\rm U}({\mathcal O}(X)))\cap {Z}^1_{{\mathcal O}(X)} 
          (G, {\rm U}({\mathcal O}(X)))_e \right)^{\perp_G}$
          for a connected algebraic group  $G$
          (for $\perp_G$, see the paragraph preceding  Lemma \ref{pairing}). 
          In the case where
           ${\mathcal O}(X)^{\Delta_{(X, G)}}$
            is finitely generated, the natural action 
            $(X\dslash \Delta_{(X, G)}, G)$ is stable.

          \begin{pr}\label{finite generation of reduced class groups}
          Suppose that ${\rm char}(K) = p = 0$ and $G$ is a  connected algebraic group such that 
         ${\mathcal O}(X)^{R_u(G)}$ is finitely generated over $K$
         as a $K$-algebra,  where $R_u (G)$ denote the unipotent radical of $G$. 
         Then $\widetilde{\rm Cl}({\mathcal O}(X), G)$ is 
         finitely generated. 
          
          \end{pr}
          
          \proof 
           By our assumption 
            the $K$-algebra ${\mathcal O}(X)^{\Delta_{(X, G)}}$
            is finitely generated (cf. \cite{Mum}), as 
            $G/R_u (G)$ is reductive.  Since 
            $\Delta_{(X, G)}\supseteq {\mathcal K}
            (X, G)$, we see that 
            \begin{eqnarray} {\frak X}(G)^{{\mathcal O}(X)^{\Delta_{(X, G)}}}
            &=& {\frak X}(G)_{{\mathcal O}(X)^{\Delta_{(X, G)}}}
            = (\Delta_{(X, G)})^{\perp_G}
            \nonumber\\
            &=&\widetilde{Z}^1_{{\mathcal O}(X)} 
          (G, {\rm U}({\mathcal O}(X)))\cap {Z}^1_{{\mathcal O}(X)} 
          (G, {\rm U}({\mathcal O}(X)))_e\label{eqn3-2}
          \end{eqnarray}
           (cf. Lemma \ref{pairing}). 
             Since ${\mathcal O}(X)^{\Delta_{(X, G)}}$
         is generated by a finite set  of  relative invariants of 
         $G$ as a $K$-algebra, there is a finite subset $\Lambda$
         of ${\frak X}(G)^{{\mathcal O}(X)^{\Delta_{(X, G)}}}$ satisfying 
         $${\frak X}(G)^{{\mathcal O}(X)^{\Delta_{(X, G)}}} \subseteq \sum_{\lambda \in \Lambda}
         \mbox{\boldmath $Z$}_0 \cdot \lambda,$$
         which shows that the condition (\ref{semigroup condition}) holds
         for ${\mathcal O}(X) = R$.   By the proof of Theorem \ref{finiteness of class group}
         the reduced class group $\widetilde{\rm Cl}({\mathcal O}(X), G)$
         is finitely generated. 
         \qqed

         \begin{th}\label{main theorem for algebraic groups}
           Under the same circumstances as in Proposition 
         \ref{finite generation of reduced class groups},    
          the following four conditions are equivalent:
          \begin{itemize}
          \item[{\it (i)}]  
          The action
          $(X\dslash {\rm tor}(v, G, \Delta_{(X, G)}), 
          G/{\rm tor}(v, G, \Delta_{(X, G)}))$ is isobaric cofree
          for some $v\in \mbox{\boldmath $N$}$. 
          \item[{\it (ii)}]  $\widetilde{{\rm Cl}}({\mathcal O}(X), G)$ is a
          finite group.
          \item[{\it (iii)}] 
          ${\rm exp} (\widetilde{\rm Cl} ({\mathcal O}(X), G)) < \infty$. 
          \item[{\it (iv)}] 
          $v := {\rm exp} (\widetilde{\rm Cl} ({\mathcal O}(X), G)) < \infty$ and
          $(X\dslash {\rm tor}(v, G, \Delta_{(X, G)}), 
          G/{\rm tor}(v, G, \Delta_{(X, G)}))$ is isobaric cofree.
          \end{itemize}
          In the case where  ${\rm U}({\mathcal O}(X)) = K^\ast$, 
          these conditions are equivalent to
          \begin{itemize}
          \item[{\it (v)}] ${\rm exp} ({\rm UrCl} ({\mathcal O}(X), G))
          = {\rm exp} (\widetilde{\rm Cl} ({\mathcal O}(X), G)) < \infty$. 
          \end{itemize}
         Especially if both $X$ and $(X, G)$ are conical, these conditions  are
         equivalent to 
         \begin{itemize}
         \item[{\it (vi)}] The natural action $(X\dslash \Delta_{(X, G)}, G)$
         is equidimensional.
         \end{itemize}  
         \end{th} 
         
         \proof
              The equivalence {\it (ii) $\Leftrightarrow$ (iii)} follows from Proposition
              \ref{finite generation of reduced class groups}. For a
              natural number $v$, by Lemma \ref{pairing} we have
              \begin{eqnarray}
              {\mathcal O}(X\dslash {\rm tor}(v, G, \Delta_{(X, G)})) &=& \bigoplus_{\chi
              \in {\rm tor}(v, G, \Delta_{(X, G)})^{\perp_G}}{\mathcal O}
              (X)_\chi\nonumber\\
              &=& \bigoplus_{\chi
              \in v\cdot {\frak X}(G)_{{\mathcal O}(X\dslash {\Delta_{(X, G)}})}}{\mathcal O}
              (X)_\chi \nonumber \\ &=& 
              \bigoplus_{\psi
              \in (\Delta_{(X, G)})^{\perp_G}}{\mathcal O}
              (X)_{v \psi}.\label{eqn3-3}\end{eqnarray}
              
              {\it (iii) $\Rightarrow$ (i) :} Put $v = 
              {\rm exp}(\widetilde{\rm Cl} ({\mathcal O}(X), G))$ and let
              $\chi$ be a cocycle in ${\rm tor}(v, G, \Delta_{(X, G)})^{\perp_G}$. Then 
              there is a cocycle $\psi\in \widetilde{Z}^1_{{\mathcal O}(X)} 
          (G, {\rm U}({\mathcal O}(X)))\cap {Z}^1_{{\mathcal O}(X)} 
          (G, {\rm U}({\mathcal O}(X)))_e$ satisfying 
          $v \psi = \chi$ (cf. (\ref{eqn3-2}) and Lemma \ref{pairing}). By the definition of
          $\widetilde{\rm Cl}({\mathcal O}(X), G)$, the divisor class
          $v\cdot [{\mathcal O}(X)_\psi]\in {\rm Cl}({\mathcal O}(X)^G)$
          equals to zero. Applying Corollary \ref{main corollary} to $\chi$, we see
          that ${\mathcal O}(X)_\chi \cong {\mathcal O}(X)^G$ as
          ${\mathcal O}(X)^G$-modules.  Thus $(X\dslash {\rm tor}(v, G, \Delta_{(X, G)}),
          G/{\rm tor}(v, G, \Delta_{(X, G)}))$ is isobaric cofree (cf. (\ref{eqn3-3})).
         The implication  {\it (iii) $\Rightarrow$ (iv)} follows from this. 
          
          {\it (i) $\Rightarrow$ (iii) :}  Let $\psi\in 
          \widetilde{Z}^1_{{\mathcal O}(X)} 
          (G, {\rm U}({\mathcal O}(X)))\cap {Z}^1_{{\mathcal O}(X)} 
          (G, {\rm U}({\mathcal O}(X)))_e$  be any cocycle. By  Lemma \ref{pairing} we have
          $v \psi \in {\rm tor}(v, G, \Delta_{(X, G)})^{\perp_G}$
          and see that ${\mathcal O}(X)_{v \psi}\cong {\mathcal O}(X)^G$
          (cf. {\it (i)}). 
          Then it follows from Corollary \ref{main corollary} that
          $v\cdot [{\mathcal O}(X)_\psi] = 0\in {\rm Cl}({\mathcal O}(X))$.
          Thus $\widetilde{\rm Cl}({\mathcal O}(X), G)$ is a torsion group, 
          which shows {\it (iii)} (cf. Proposition 
          \ref{finite generation of reduced class groups}). 
          
          Suppose  ${\mathcal O}(X)$ has only trivial units. Then
          $Z^1(G, {\rm U}({\mathcal O}(X))) \cong H^1(G, {\rm U}({\mathcal O}(X)))$.  
          As $H^1(G, {\rm U}({\mathcal O}(X)))_{R, 0}
          \subseteq {\frak X}(G) \subseteq H^1(G, {\rm U}({\mathcal O}(X)))$, 
          by Remark \ref{zero cohomology} the equivalence {\it (ii) $\Leftrightarrow$
          (v)} follows  from  Theorem \ref{finiteness of class group}.
           
           Suppose that both $X$ and $(X, G)$ are conical. The implication
           {\it (i) $\Rightarrow$ (vi)} is a consequence of the generic
           fiber theorem 
           of graded version of flat local morphisms (e.g., \cite{Matsumura}). 
           So suppose that the condition {\it (vi)} holds. 
           Let  $$\psi\in 
          \widetilde{Z}^1_{{\mathcal O}(X)} 
          (G, {\rm U}({\mathcal O}(X)))\cap {Z}^1_{{\mathcal O}(X)} 
          (G, {\rm U}({\mathcal O}(X)))_e= (\Delta_{(X, G)})^{\perp_G}$$ 
          be a cocycle. Since  
           $(X\dslash \Delta_{(X, G)}, G/\Delta_{(X, G)})$
          is a stable and equidimensional action of an algebraic torus
          on a conical normal variety,   applying Sect. 4 
         of \cite{Nak1} to ${\mathcal O}(X)_\psi$, we can choose  
         a natural number $r$ in such a way that  
         $${\mathcal O}(X)_{n r \psi }= 
         {\mathcal O}(X\dslash \Delta_{(X, G)})_{n r \psi}
         \cong {\mathcal O}(X)^G$$ for all $n \in \mbox{\boldmath $N$}$
         as ${\mathcal O}(X)^G$-modules. Hence $\widetilde{{\rm Cl}}({\mathcal O}
         (X), G)$ is a torsion group (cf. Corollary \ref{main corollary}), which shows
         {\it (ii)}. 
          \qqed

          \begin{cor}\label{reduction case}   Under the same circumstances as in Proposition 
         \ref{finite generation of reduced class groups}, suppose that
         both $X$ and $(X, G)$ are conical.  Then the 
         following conditions are equivalent:
         \begin{itemize}
         \item[\it (i)] $v:={\rm exp} ({\rm UrCl} ({\mathcal O}(X\dslash {\mathcal K}(X, G)), G))
         < \infty$ and $(X\dslash {\rm tor}(v, G, {\frak  R}(X\dslash {\mathcal K}
         (X, G), G)), G)$
         is cofree. 
         \item[\it (ii)] ${\rm exp} ({\rm UrCl} ({\mathcal O}(X\dslash {\mathcal K}(X, G)), G))
          = {\rm exp} (\widetilde{\rm Cl} ({\mathcal O}(X\dslash {\mathcal K}(X, G)), G)) < \infty$
          and the inclusion ${\mathcal O}(X)^G \to {\mathcal O}(X)^{{\mathcal K}(X, G)}$
          is no-blowing-up of codimension one (cf. p. 30 of \cite{Fossum}).
          \item[\it (iii)] $\widetilde{\rm Cl} ({\mathcal O}(X\dslash {\mathcal K}(X, G)), G)$
          is a finite group and the inclusion 
          ${\mathcal O}(X)^G \to {\mathcal O}(X)^{{\mathcal K}(X, G)}$ 
          is no-blowing-up of codimension one. 
         \item[{\it (iv)}] $(X\dslash {\mathcal K}(X, G), G)$
         is equidimensional. 
        \end{itemize}  
          \end{cor}
          
          \proof
          In order to show this corollary, we may suppose that
          ${\mathcal K}(X, G) =\{1\}$ as $${\rm tor}(m, G, {\frak  R}(X\dslash {\mathcal K}
         (X, G), G))\supseteq {\mathcal K}(X, G)$$ 
          for $m \in\mbox{\boldmath $N$}$. So
           $(X, G)$ is a faithful stable action of an algebraic torus $G$. The
           implication {\it (i)} $\Rightarrow$ {\it (iv)} and 
           the equivalence {\it (iv)} $\Leftrightarrow$ ``$(X\dslash {\frak R}(X, G),
           G)$ is equidimensional" follow from the finiteness of 
           ${\frak R}(X, G)$.  Each  condition of  {\it (i)} - {\it (iv)} 
           satisfies  the condition that
            the inclusion ${\mathcal O}(X)^G \to {\mathcal O}(X)^{{\mathcal K}(X, G)}$
          is no-blowing-up of codimension one 
          (cf. (PDE), p. 30 of \cite{Fossum}), which is assumed in the following proof.
         Then  $\widetilde{Z}^1_{{\mathcal O}(X)}(G, {\rm U}({\mathcal O}(X))
         ={\frak X}(G)_{{\mathcal O}(X)} = {\frak X}(G)$. 
             On the other hand,  by Proposition 
           \ref{pseudo-reflection group},  ${\frak R}(X, G)^{\perp_G} =
           Z^1_{{\mathcal O}(X)}(G, {\rm U}({\mathcal O}(X))_e$, which shows
           \begin{eqnarray} \Delta_{(X, G)} = 
           \left(\widetilde{Z}^1_{{\mathcal O}(X)}(G, {\rm U}({\mathcal O}(X))
           \cap Z^1_{{\mathcal O}(X)}(G, {\rm U}({\mathcal O}(X))_e
           \right)^{\perp_G} = {\frak R}(X, G)\nonumber
           \end{eqnarray}
           (cf. Lemma \ref{pairing}). Thus the equivalence  
           of conditions in this corollary
           is a consequence of Theorem \ref{main theorem for algebraic groups}. \qqed
           
           \medskip
           
          Suppose that  $G$ is an
          algebraic torus and 
          $(X, G)$ is stable.  Clearly ${\mathcal Q}({\mathcal O}(X)^G )
          = {\mathcal Q}({\mathcal O}(X))^G $. 
          Moreover, in the
           case where  $t_{(X,G)} :=
           {\rm exp} ({\rm UrCl}({\cal O}(X), G))$ is finite, express
           \begin{equation} t_{(X, G)} = \tilde{t}_{(X, G)} \cdot t^{{\frak R}}_{(X, G)}
           \label{3.1}
           \end{equation}
           as a product of natural numbers $\tilde{t}_{(X, G)}$,
            $t^{{\frak R}}_{(X, G)}$ 
            satisfying that 
           $${\rm GCD}(\tilde{t}_{(X, G)},
           \left\vert {\frak R}(X, G)\vert_X\right\vert) = 1$$ and 
           $ \left\vert {\frak R}(X, G)\vert_X\right\vert$ is divisible in $\mbox
           {\boldmath $Z$}$ by
           any prime divisor of $t^{{\frak R}}_{(X, G)}$.  
         Let $F_{(X, G)}$ be the subgroup of ${\frak R}(X, G)$
         consisting all elements $\sigma\in {\frak R}(X, G)$
          such that prime divisors of 
         ${\rm ord}(\sigma\vert_X)$ are divisors of $t^{\frak R}_{(X, G)}$. 
         Obviously 
           there exists a natural number $k$ such that 
           \begin{eqnarray}
           {\rm tor} ((t^{{\frak R}}_{(X, G)})^k, G, L_{(X, G)}) \cap {\frak R}(X, G)
          &=& {\rm tor} ((t^{{\frak R}}_{(X, G)})^{k+j}, G, L_{(X, G)}) \cap {\frak R}(X, G)
          \nonumber\\
          &=& F_{(X, G)}\nonumber\end{eqnarray}
           for any $j \in \mbox{\boldmath $N$}$.

            \begin{df}\label{obstruction}
              \rm Under the same  circumstances as above, define
            the {\it obstruction subgroup for cofreeness of}  $(X, G)$, denoted
           by ${\rm Obs}(X, G)$, 
           as follows; 
            $${\rm Obs}(X, G) := \widetilde{{\frak  R}}(X\dslash H_{(X, G)} 
           , G)$$ 
           where $H_{(X, G)} := {\rm tor}(\tilde{t}_{(X, G)}, G, L_{(X, G)})\cdot 
          {\rm tor}(t^{\frak R}_{(X, G)}, G, F_{(X, G)})$ (cf. Definition 
          \ref{def of pseudo-reflection group}). 
            \end{df} 
          \begin{rem} \rm ~Obviously     ${\rm Obs}(X, G)\vert_X$ is  a finite
           group. If $X$ (i.e.,  ${\mathcal O}(X)$)
          is factorial, then $t_{(X, G)} =1$ and we see that $L_{(X, G)}    
          = {\rm Obs}(X, G)$, i.e., ${\rm Obs}(X, G)\vert_X = \{1\}$. 
          On the other hand in the case where
          ${\rm U}({\mathcal O}(X)) = K^\ast$,  unless ${\rm Obs}(X, G)\vert_X = \{1\}$, 
          $(X, G)$ is never isobaric cofree (cf. Corollary
          \ref{main corollary} and Theorem \ref{main theorem for algebraic groups}). 
          \end{rem}
          
              {\it There are many examples of equidimensional actions of (connected) algebraic  tori on conical normal varieties
          which are not cofree} as follows:

          \begin{exam}\label{example1}\rm
           Let $V = K^4$ 
           be a $4$-dimensional vector
          space over $K$ and let $M$ be a subgroup
          \[\left<\left\{\left.\left(\begin{array}{cccc}
          u & & & \\
          & u^{-1} & &\\
          & &  v & \\
          & & & v^{-1}
          \end{array}\right)~\right\vert~ (u, v)\in (K^\ast)^2 \right\}
          \bigcup \left\{\tau =\left(\begin{array}{cccc}
          \zeta_3 & & & \\
          & \zeta_3 & &\\
          & &  \zeta_3^{-1} & \\
          & & & \zeta_3^{-1}
          \end{array}\right)\right\}\right> \]
          of $GL(V^\vee)$ 
          where $V^\vee$ the $K$-dual space of $V$  and $\zeta_3$ is a 
          primitive $3$-th root of $1$ in $K$. Put $X:= V\dslash \langle \tau\rangle$
          and $G := M^0$ which acts naturally on $X$. Since $M\subseteq SL(V)$,
          we see that ${\frak R}(V, M) =\{1\}$ and hence ${\frak R}(X, G) = \{1\}$. 
          Clearly $V\dslash G \cong \mbox{\boldmath $A$}_K^2$ and 
          the action $(V, G)$ is cofree, which implies that
          $(X, G)$ is equidimensional. Applying
          Samuel's Galois descent to the action $(V\dslash G, \langle \tau \rangle)$,
          we have $${\rm Cl}({\mathcal O}(X)^G) ={\rm Cl}(({\mathcal O}(V)^G)^{\langle
          \tau\rangle})  \cong \mbox{\boldmath $Z$}/
          3\mbox{\boldmath $Z$}.$$ Because $(X, G)$ is not cofree,
          we see $\{[\theta]\}\not= \widetilde{\rm Cl} ({\mathcal O}(X), G)
          = {\rm Cl}({\mathcal O}(X)^G)$ and 
          $t_{(X, G)} = 3$ (recall $\theta$
          denotes the trivial cocycle).  
          Consequently $${\rm Obs}(X, G)\vert_X \cong \mbox{\boldmath $Z$}/
          3\mbox{\boldmath $Z$} \oplus \mbox{\boldmath $Z$}/
          3\mbox{\boldmath $Z$}.$$ 
          
          \end{exam}
          
          \begin{exam}\label{example2} \rm Let $V = K^4$ and put
          \[G:= \left\{\left.\left(\begin{array}{cccc}
          tu & & & \\
          & tu^{-1} & &\\
          & &  t & \\
          & & & t^{-3}
          \end{array}\right)~\right\vert~ (t, u)\in (K^\ast)^2 \right\} \subseteq GL(V^\vee),\]
          \[H:= \left\{\left.\left(\begin{array}{cccc}
          t & & & \\
          & t & &\\
          & &  t & \\
          & & & t^{-3}
          \end{array}\right)~\right\vert~ t\in K^\ast \right\} \subseteq GL(V^\vee),\]
          and $X : = V\dslash H$.  Here the matrix representation
          is given by the basis $\{X_1, \dots, X_4\}$. 
          As in Example \ref{example1}, we similarly see that ${\frak R}(X, G) = \{1\}$. 
          By \cite{Nak4} one sees that ${\rm Cl}({\mathcal O}(X))
          \cong H_{X_4} \cong \mbox{\boldmath $Z$}/3\mbox{\boldmath $Z$}$.
          There is a finite dominant morphism $\pi : X \to Y =\mbox{\boldmath $A$}^3$, 
          where $Y$ is defined by 
          $K[X_1^3X_4, X_2^3X_4, X_3^3X_4]$ and  $\pi$ is associated with 
          $ K[X_1^3X_4, X_2^3X_4, X_3^3X_4]\hookrightarrow {\mathcal O}(X)$. 
          Clearly $G$ acts naturally on $Y$ and $\pi$ is $G$-equivariant. 
          Since $(Y,G)$ is cofree, $(X, G)$ is stable and equidimensional. On the other hand
          we easily see that both $(V, G)$ and
           $(X, G)$ are not cofree. Hence $\widetilde{\rm Cl}
          ({\mathcal O}(X), G) \not=\{0\}$ and 
          ${\rm Cl}({\mathcal O}(X))
          \supseteq {\rm UrCl}({\mathcal O}(X), G) 
          \cong \mbox{\boldmath $Z$}/3\mbox{\boldmath $Z$}$.
           Consequently $${\rm Obs}(X, G)\vert_X \cong \mbox{\boldmath $Z$}/
          3\mbox{\boldmath $Z$}.$$ 
          Obviously  $X$ is regarded as an affine toric variety with a non-cofree
          equidimensional torus action commuting with its toric structure. 
          In general, for any finitely generated abelian group $A$, 
          there exists an affine toric variety whose class group
          is isomorphic to $A$ (cf. \cite{Nak4}). 
          \end{exam}

         \begin{th}\label{equi}
           Suppose that $G$ is a connected algebraic torus
           and $p = 0$. Suppose that both $X$ and $(X, G)$ are conical. If $(X, G)$
           is stable, then the following 
         conditions are equivalent:

         \begin{itemize}
         \item[{\it (i)}] The action $(X, G)$ is equidimensional.
         \item[{\it (ii)}] The exponent ${\rm exp} ({\rm UrCl}({\mathcal O}(X), G))$  
          is finite and
         the    action $(X\dslash {\rm Obs}(X, G), G) =
         (X\dslash {\rm Obs}(X, G), G/{\rm Obs}(X, G))$
          is cofree.

         \end{itemize}
         
         \end{th}

         \proof {\it (ii)} $\Rightarrow$ {\it (i) :}~ Since cofree conical actions are
         equidimensional and the group ${\rm Obs}(R, G)\vert _X$ is finite,
         we immediately see that {\it (ii)} implies {\it (i)}. 
         
         {\it (i)} $\Rightarrow$ {\it (ii) :}~ Suppose that $(X, G)$
         is equidimensional. We use notation in  Definition \ref{obstruction}
         and the paragraph preceding to
         Definition \ref{obstruction}. By  Proposition \ref{generalized reflection}
         we see that $${\rm Obs}(X, G) \lhd {\frak R}(X\dslash H_{(X, G)}, G) = 
         H_{(X, G)} \cdot {\frak R}(X, G)$$ and $(X\dslash {\rm Obs}(X, G), 
         {\frak R}(X\dslash H_{(X, G)}, G))$  is cofree. Let $\chi \in {\frak X}(G)$ 
          be a character satisfying $\chi (H_{(X, G)} \cdot {\frak R}(X, G))
         = \{1\}$. Note $\chi \in {\frak X}(G)_{{\mathcal O}(X)}$, as $\chi (L_{(X, G)}) =\{1\}$.  
         It suffices to show that ${\mathcal O}(X)_\chi \cong {\mathcal O}(X)^G$
         as ${\mathcal O}(X)^G$-modules. 
             Let $\sigma$ be any element of $G$
         such that $\sigma^ {t_{(X, G)}} \in {\frak R}(X, G)$. 
         Recall the expression (\ref{3.1}). 
         As the map
         $${\frak R}(X, G)\ni \tau \to \tau^{\tilde{t}_{(X, G)}} \in {\frak R}(X, G)$$
         induces an automorphism of ${\frak R}(X, G)\vert_X$, the 
         element $\sigma^{t^{\frak R}_{(X, G)}}$ belongs to the subgroup 
         ${\rm tor}(\tilde{t}_{(X, G)}, G, L_{(X, G)}) \cdot {\frak R}(X, G)$. 
            Since
         $${\rm tor}(\tilde{t}_{(X, G)}, G, L_{(X, G)})\ni \tau \to \tau^{t^{\frak R}_{(X, G)}}
         \in {\rm tor}(\tilde{t}_{(X, G)}, G, L_{(X, G)})$$ induces 
         an automorphism of ${\rm tor}(\tilde{t}_{(X, G)}, G, L_{(X, G)})\vert_X$
         and ${\frak R}(X, G)$ contains  $L_{(X, G)}$, we can choose $\mu$ from
         ${\rm tor}(\tilde{t}_{(X, G)}, G, L_{(X, G)})$ in such a way that 
         $(\sigma\cdot \mu)^{t^{\frak R}_{(X, G)}}\in {\frak R}(X, G)$. Then express 
         $(\sigma\cdot \mu)^{t^{\frak R}_{(X, G)}} = \gamma_1\cdot \gamma_2$
         for some $\gamma_i\in {\frak R}(X, G)$ such that
         $${\rm GCD} ({\rm ord} (\gamma_1\vert_X), t^{\frak R}_{(X, G)}) = 1$$
         and any prime divisor of ${\rm ord} (\gamma_2\vert_X)$ is a 
         divisor of $\left\vert {\frak R}(X, G)\vert _X\right\vert$.
         Clearly   $ \gamma_2\in F_{(X, G)}$. 
         Since $F_{(X, G)} \supseteq L_{(X, G)}$
         and $\langle \gamma_1\rangle \ni \tau \to \tau^{t^{\frak R}_{(X, G)}}
         \in \langle \gamma_1\rangle$ induces an automorphism of 
         $\langle \gamma_1\rangle\vert_X$, there exists an element $\delta\in 
         \langle \gamma_1\rangle$ satisfying 
         $$(\sigma \cdot \mu\cdot \delta)^{t^{\frak R}_{(X, G)}}\in F_{(X, G)}.$$
         Consequently we see $$\sigma \in {\rm tor}(\tilde{t}_{(X, G)}, G, L_{(X, G)})
         \cdot {\rm tor}(t^{\frak R}_{(X, G)}, G, F_{(X, G)})\cdot {\frak R}(X, G)
         = H_{(X, G)} \cdot {\frak R}(X, G).$$
         Thus $\chi (\sigma) =1$, which implies $$\chi ({\rm tor}(t_{(X,G)},
         G, {\frak R}(X, G))) = \{1\}.$$ By 
         the implication
         {\it (iv)} $\Rightarrow$ {\it (i)} of 
        Corollary  \ref{reduction case}, we see
         that ${\mathcal O}(X)_\chi \cong {\mathcal O}(X)^G$ as 
         ${\mathcal O}(X)^G$-modules.  
         \qqed 
         
         \medskip

         \noindent{\it Proof of Theorem 1.2.}~ Suppose
         that $(X, G)$ is faithful and by Theorem \ref{equi} we may suppose that 
         $$t_{(X, G)}:={\rm exp}\left({\rm UrCl}({\mathcal O}(X), G)\right)
         = {\rm exp}\left(\widetilde{{\rm Cl}}({\mathcal O}(X), G)\right) < \infty.$$
          By Definition \ref{obstruction}
         we see that prime divisors of 
         $\vert{\rm Obs}(X, G)\vert$ are same as those of  $t_{(X, G)}$, which implies that 
          there exists a power of  $t_{(X, G)}$ which 
         is divisible by $\vert{\rm Obs}(X, G)\vert$. 
         \qqed
         
          \medskip
          
         Now, Corollary 1.3 follows immediately from the definition of 
         $\widetilde{{\rm Cl}}({\mathcal
         O}(X), G)$ and Theorem \ref{equi}. 
         
         \vskip 0.4cm


        \end{document}